\documentclass[11pt]{amsart}
\usepackage{latexsym,url}
\usepackage[all]{xy}

\numberwithin{equation}{section}

\newtheorem{theo}{Theorem}
\swapnumbers

\theoremstyle{plain}
\newtheorem{prp}[subsection]{Proposition}
\newtheorem{lem}[subsection]{Lemma}
\newtheorem{lemu}[subsubsection]{Lemma} 

\theoremstyle{definition}

\newtheorem{rems}[subsection]{Remarks}
\newtheorem{remu}[subsubsection]{Remark}

\newcommand{\bsm}{\begin{smallmatrix}}
\newcommand{\esm}{\end{smallmatrix}}

\newenvironment{pf}{\medskip\noindent{\em Proof:}}{\hspace*{\fill}$\Box$}

\newcommand{\df}{\colon}
\newcommand{\ra}{\rightarrow}             
\newcommand{\id}{1\kern -.35em 1}
\newcommand{\ebrace}[1]{\langle #1\rangle}
\newcommand{\abs}[1]{\lvert #1\rvert}
\newcommand{\dimv}{\operatorname{\underline{dim}}}

\newcommand{\lam}{\lambda}
\newcommand{\iot}{\iota}
\newcommand{\alp}{\alpha}
\newcommand{\vph}{\varphi}
\newcommand{\vpi}{\varpi}
\newcommand{\ome}{\omega}

\newcommand{\Hom}{\operatorname{Hom}}
\newcommand{\Ext}{\operatorname{Ext}}
\newcommand{\End}{\operatorname{End}}
\newcommand{\GL}{\operatorname{GL}}
\newcommand{\SL}{\operatorname{SL}}
\newcommand{\Obj}{\operatorname{Obj}}
\newcommand{\Lie}{\operatorname{Lie}}  

\newcommand{\ka}{\mathrm{k}}        
\newcommand{\EN}{\mathbb{N}}
\newcommand{\ENn}{\EN_0}            
\newcommand{\Zet}{\mathbb{Z}}
\newcommand{\CC}{\mathbb{C}}        

\newcommand{\Aa}{\mathsf{A}}             
\newcommand{\Bb}{\mathsf{B}}
\newcommand{\Cc}{\mathsf{C}}
\newcommand{\Dd}{\mathsf{D}}
\newcommand{\Ee}{\mathsf{E}} 
\newcommand{\Ff}{\mathsf{F}}
\newcommand{\Gg}{\mathsf{G}}

\newcommand{\ZQ}{\Zet Q}

\newcommand{\AR}[1]{A_{#1}}         
\newcommand{\cA}[1]{\check{A}_{#1}}
\newcommand{\tAi}{\widetilde{A}_\ii}
\newcommand{\bQ}{\overline{Q}}      

\newcommand{\de}{\mathbf{e}}
\newcommand{\dv}{\mathbf{v}}             
\newcommand{\dw}{\mathbf{w}}
\newcommand{\di}{\mathbf{i}}
\newcommand{\dpr}{\mathbf{p}}

\newcommand{\lmd}[1]{{#1}\operatorname{-mod}}
\newcommand{\lind}[1]{{#1}\operatorname{-ind}}
\newcommand{\slmd}[1]{{#1}\text{-}\underline{\operatorname{mod}}}
\newcommand{\srmd}[1]{\underline{\operatorname{mod}}\text{-}#1}

\newcommand{\linj}[1]{{#1}\operatorname{-inj}}
\newcommand{\cD}{\mathcal{D}}          
\newcommand{\Der}[1]{\mathcal{D}^b(#1)}

\newcommand{\Lam}{\Lambda}          
\newcommand{\tLam}{\widetilde{\Lam}}
\newcommand{\kQ}{{\ka Q}}           
\newcommand{\Gam}{\Gamma}
\newcommand{\AuC}[1]{\Gam_{#1}}     
\newcommand{\cGam}{\check{\Gam}}    
\newcommand{\uGam}{\underline{\Gam}}
\newcommand{\huGam}{\widehat{\uGam}}

\newcommand{\g}{\mathfrak g}    
\newcommand{\h}{\mathfrak h}
\newcommand{\n}{\mathfrak n}

\newcommand{\Nak}{\mathcal{N}}     
\newcommand{\Sig}{\Sigma}          
\newcommand{\Ome}{\Omega}          
\newcommand{\hnu}{\hat{\nu}}       

\newcommand{\Jx}{J^{\textstyle{\cdot}}}
\newcommand{\Jr}{J^{\rho}}
\newcommand{\RA}{\operatorname{R}} 
\newcommand{\sM}{\operatorname{I}} 
\newcommand{\sP}{\operatorname{P}} 

\newcommand{\M}{\mathcal{M}}        
\newcommand{\cE}{\mathcal{E}}       
\newcommand{\ii}{\mathbf{i}}        
\newcommand{\f}{\mathfrak f}        
\newcommand{\Del}{\Delta}           
\newcommand{\pRt}{\Pi}              
\newcommand{\bcsl}{\backslash}
\newcommand{\op}{^{\text{op}}}

\title{Auslander algebras and initial  seeds for cluster algebras}

\author{Christof Gei{\ss}}
\address{Instituto de Matem\'aticas, UNAM, Ciudad Universitaria,
04510 M\'exico D.F., Mexico}
\email{christof@math.unam.mx}

\author{Bernard Leclerc}
\address{Laboratoire LMNO, Universit\'e de Caen,
F-14032 Caen Cedex, France}
\email{leclerc@math.unicaen.fr}

\author{Jan Schr\"oer}
\address{Department of Pure Mathematics, University of Leeds,
Leeds LS2 9JT, UK}
\email{jschroer@maths.leeds.ac.uk}

\date{August 11, 2006}
\thanks{Ch.G. acknowledges support from DGAPA grant IN101402-3.
B.L. is grateful to the GDR 2432 and the GDR 2249 for their support. 
J.S. was supported by a research fellowship from the DFG}
\begin{document}
\begin{abstract}
Let $Q$ be a Dynkin quiver and $\pRt$ the corresponding set of positive
roots. For the preprojective algebra $\Lam$ associated to $Q$ we produce
a rigid $\Lam$-module $\sM_Q$ with $r=\abs{\pRt}$ pairwise non-isomorphic
indecomposable direct summands by pushing the injective modules of the 
Auslander algebra of $\kQ$ to $\Lam$.

If $N$ is a maximal unipotent subgroup of a complex simply connected simple
Lie group of type $\abs{Q}$, then the coordinate ring $\CC[N]$ is an
upper cluster algebra. 
We show that the elements of the 
dual semicanonical basis which correspond to the indecomposable direct summands
of $\sM_Q$ coincide with certain generalized minors which 
form an initial cluster for $\CC[N]$, and that the corresponding exchange
matrix of this cluster can be read from the Gabriel quiver of 
$\End_\Lam(\sM_Q)$.

Finally, we exploit the fact that the categories of injective modules over
$\Lam$ and over its covering $\tLam$ are triangulated in order to show
several interesting identities in the respective stable module categories.
\end{abstract}
\maketitle

\section*{Introduction}
Let $\ka$ be a field and $Q$ be a Dynkin quiver. So the underlying 
graph $\abs{Q}$ of $Q$ is a simply laced Dynkin diagram.
We produce for the preprojective algebra $\Lam$ over $\ka$ associated to 
$Q$ a module $\sM_Q$ by pushing a minimal injective cogenerator over the
Auslander algebra $\AuC{Q}$ of $\kQ$ to $\lmd{\Lam}$. 
It is easy to see that $\sM_Q$ decomposes into $r=\abs{\Pi}$ pairwise 
non-isomorphic direct summands. We show that $\sM_Q$ is a rigid module 
module, i.e. $\Ext_\Lam^1(\sM_Q,\sM_Q)=0$. Moreover, the Gabriel  quiver 
$\cA{Q}$  of $\End_\Lam(\sM_Q)\op$ is obtained from the 
Auslander-Reiten quiver $\AR{Q}$ of $\kQ$ by inserting an extra arrow 
$x\ra\tau x$ for each non-projective
vertex $x$.

In~\cite{gesr02} we have shown that if $M=\oplus_{i=1}^m M_i$ for pairwise
non-isomorphic indecomposable $\Lam$-modules $M_i$,  then 
$\Ext_\Lam^1(M,M)=0$ implies $m\leq r$. So our result shows that
this maximum is assumed for each Dynkin quiver. 
By~\cite[Theorem 2.2]{GeLeSc05a} we conclude that $\sM_Q$ is a maximal
1-orthogonal $\Lam$-module and thus $\End_\Lam(\sM_Q)$ is
a higher Auslander algebra in the sense of Iyama~\cite{Iyam04a}. 
It follows that each rigid module $M$ as above can be completed 
to a rigid module with $r$ pairwise non-isomorphic indecomposable direct
summands.
Note, that for the proof of the main result
of~\cite{GeLeSc05a} it is essential that the quiver $\cA{Q}$ has no loops.

Let now $G$ be a complex simply connected simple  Lie group of
type $\abs{Q}$ with $N\subset G$ a maximal unipotent subgroup. 
Choose $\ii=(i_1,i_2,\ldots,i_r)$ a reduced expression for the longest element
$w_0$ of the Weyl group $W$ of $G$.
It follows from~\cite{BeFoZe05} that the coordinate
ring $\CC[N]$ is  an (upper) cluster algebra. 
Associated to $\ii$ one obtains an initial seed 
$(\Del(j,\ii)'_{j=1,2,\ldots, r},\widetilde{B}(\ii)')$,
where the $\Del(j,\ii)'$ are certain generalized minors,
and the exchange matrix $\widetilde{B}(\ii)'$ is 
obtained naturally from the quiver $\cA{Q}$ described above.

Next, $\ii$ provides us with a convenient labelling of the indecomposable 
direct  summands of $\sM_Q$, that is, $\sM_Q=\oplus_{j=1}^r\sM(j,\ii)$. 
Clearly the $\sM(j,\ii)$ are rigid. Thus, if we restrict to the special case 
$\ka=\CC$ the $\Lam$-modules $\sM(j,\ii)$ serve as natural
labels for  elements $\rho_{\sM(j,\ii)}$ of the dual of Lusztig's semicanonical
basis. This is a natural  basis for $\CC[N]$ and  we show that 
$\rho_{\sM(j,\ii)}=\Del(j,\ii)'$. 

\section{Main results}
\subsection{} \label{ssec:mr1}
We say that a quiver $Q$ is a {\em Dynkin quiver} if its
underlying graph $\abs{Q}$ is a Dynkin diagram of type $\Aa,\Dd,\Ee$.
For a field $\ka$ we consider the {\em path category} $\ka[Q]$ (or
$\kQ$ for short). The category $\lmd{\kQ}$ of finitely presented $\ka$-functors
$M\!\df\kQ\ra\lmd{\ka}$ is equivalent to the category of finitely presented
left modules over the corresponding path algebra which we denote by some
abuse also by $\kQ$.

For a quiver $Q$ we consider the double $\bQ$ which is obtained from $Q$
by adding a new arrow $i\xleftarrow{a^*} j$ for each arrow $
i\xrightarrow{a}j$ in $Q$. The preprojective  algebra $\Lam$ is
the quotient of the path algebra $\ka\bQ$ by  
the ideal generated by the elements
\[
\rho_q =  \sum_{\substack{a\in Q_1\\ t(a)=q}} a^* a
         -\sum_{\substack{a\in Q_1\\ h(a)=q}} a a^*
\qquad\text{for}\quad q\in Q_0,
\]
see also~\ref{ssec:quica}. 

In what follows, $Q$ will always be a connected Dynkin quiver.  
This implies that $\Lam$ is a finite-dimensional selfinjective algebra,
which depends only on $\abs{Q}$. Like $\kQ$, the algebra $\Lam$ can also be
considered as a $\ka$-category.

We have the universal covering $F\df\tLam\ra\Lam$ where $\tLam$ is the
path category $\ka[\ZQ]$ modulo the usual mesh relations. 
The fundamental group $\Zet$ of $\Lam$ acts on $\tLam$  via the 
translation $\tau$.
Associated to $F$ we have the push-down functor
$F_\lam\df\lmd{\tLam}\ra\lmd{\Lam}$, see~\ref{ssec:cov-constr}.

In $\ZQ$ we find the Auslander-Reiten quiver $\AR{Q}$ of $\kQ$ as a full
convex subquiver. The Auslander category $\AuC{Q}$ is the full subcategory of 
$\tLam$ which has the vertices of $\AR{Q}$ as objects. Denote the inclusion
of $\AuC{Q}$ into $\tLam$  by $J$. 
There is a natural equivalence $\RA_Q$ from $\AuC{Q}$ to the category of
indecomposable $\kQ$-modules, $\lind{\kQ}$.
We say that an object $x\in\Obj(\AuC{Q})$ is projective if 
$\tau x\not\in\Obj(\AuC{Q})$, dually $x$ is injective if 
$\tau^{-1}x\not\in\Obj(\AuC{Q})$,
see~\ref{ssec:ausca}.

Associated to $\AuC{Q}$ we consider the $\ENn$-graded category 
$\cGam_Q$. It has the same objects as $\AuC{Q}$ but the morphisms of
degree $i$ are given by $\cGam_{Q,i}(x,y) =\AuC{Q}(\tau^i x, y)$ 
if $\tau^i x\in\Obj(\AuC{Q})$. We equip $\cGam_Q$
with the natural composition. The Gabriel quiver $\cA{Q}$ of $\cGam_Q$ is
obtained from $\AR{Q}$ by inserting an additional (degree $1$)
arrow $x\ra \tau x$ for each non-projective $x\in\Obj(\AuC{Q})$,
see~\ref{ssec:start-grad}.

\subsection{Start modules.} \label{ssec:startM}
Let us write  $D$ for the usual duality $\Hom_\ka(-,\ka)$.
We denote by $\Jx$ the functor which considers a $\AuC{Q}$-module (trivially) 
as a $\tLam$-module. Thus if apply the functor $F_\lam\Jx$ to
the  injective $\AuC{Q}$-module $D\AuC{Q}(-,x)$ we obtain
a $\Lam$-module. We call
\[
\sM_Q:= \bigoplus_{x\in\Obj(\AuC{Q})} F_\lam \Jx (D\AuC{Q}(-,x)).
\]
the {\em start module} for $\Lam$ associated to $Q$. Note that
$F_\lam \Jx (D\AuC{Q}(-,  x))$ is isomorphic to a submodule of
$F_\lam \Jx (D\AuC{Q}(-, \tau^{-1} x))$ if $x\in\Obj(\AuC{Q})$ is not 
injective,
and $F_\lam \Jx (D\AuC{Q}(-, x))$ is an injective $\Lam$-module
if $x\in\Obj(\AuC{Q})$ is injective, see~\ref{ssec:ausca} and~\ref{ssec:adj}.

Consider $\cE=\oplus_{x,y\in{\Obj(\AuC{Q})}}\cGam_Q(x,y)$ as a
(graded) associative $\ka$-algebra  with multiplication induced from the
composition of morphisms. 

\begin{theo} \label{thm:start}
Let $\Lam$ be the preprojective algebra associated to a Dynkin quiver $Q$.
Then $\mathcal{E}$ is isomorphic to
$\End_\Lam(\sM_Q)\op$. In particular,
the Gabriel quiver of $\End_\Lam(\sM_Q)\op$ is identified with $\cA{Q}$ 
as described above in~\ref{ssec:mr1}. Moreover $\sM_Q$ is rigid in
the sense that
$\Ext^1_\Lam(\sM_Q,\sM_Q)=0$.
\end{theo}

The proof of this result is prepared in~\ref{ssec:start-prep}, 
\ref{prp:form}, \ref{prp:start-eqn} and finished
in~\ref{ssec:start-pf}.

\subsection{Reduced expressions.} \label{intro:expr}
Let $\pi\df \Obj(\AuC{Q})\ra Q_0$ denote the map induced by the composition
$FJ$.
We call a total ordering $x(1)<x(2)<\ldots <x(r)$ of the objects of
$\AuC{Q}$ {\em adapted (to $Q$)} if $\AuC{Q}(x(i),x(j))=0$ for $i<j$. 
It is easy to find such orderings given that the quiver $A_Q$ is directed.

We call a vertex $i\in Q_0$ a {\em source} in $Q$ if no arrow ends at $i$.
In this case we denote by $s_i(Q)$ the quiver which is obtained from $Q$ by
reversing each arrow starting in $i$.

Let $\ii=(i_1,i_2,\ldots,i_r)$ be a reduced expression for the longest element
$w_0\in W$, that is $w_0=s_{i_1} s_{i_2}\cdots s_{i_r}$ where
$r=\abs{\pRt}$. We say that
$\ii$ is  {\em adapted to} $Q$ if $i_1$ is a source in $Q$ and
$i_{k+1}$ is a source in $s_{i_{k}}\cdots s_{i_1}(Q)$ if $1\leq k<r$.
This is  dual to the original definition in~\cite{lusz90}.

If $x(1)<\cdots <x(r)$ is an adapted ordering then
$\ii:=(\pi(x(1)),\ldots,\pi(x(r)))$ is a reduced 
expression for the longest element $w_0$ of the Weyl group $W$ associated
to $\abs{Q}$, see~\ref{intro:Liea} below. In fact,
the adapted orderings of $\Obj(\AuC{Q})$ correspond in this way bijectively 
to the reduced  expressions for $w_0$  which are  adapted to $Q$, 
see~\cite[Theorem 2.5]{beda00}.
We  set
\[
\sM(j,\ii):= F_\lam\Jx D\AuC{Q}(-,x(j))\quad\text{ for } 1\leq j\leq r.
\]
Thus, $\ii$ provides us with a convenient way of labelling the direct summands
of $\sM_Q$. 

\subsection{} \label{intro:Liea}
Let now $\g$ be a complex simple Lie algebra of type $\abs{Q}$
with the usual Serre generators $e_i,h_i,f_i$ for $i\in Q_0$. Thus the
$h_i$ form a basis of the abelian subalgebra $\h$ and the $e_i$ resp.~$f_i$
generate maximal nilpotent subalgebras $\n$ resp.~$\n_-$. The simple roots
$\alp_i$ form a basis of the dual space $\h^*$ such that $\alp_i(h_j)=a_{i,j}$
where $(a_{i,j})_{i,j\in {Q_0}}$ is the Cartan matrix of $\abs{Q}$. The
fundamental weights $(\vpi_i)_{i\in {Q_0}}$ are the basis of $\h^*$ dual
to the basis $(h_i)_{i\in Q_0}$ of $\h$.

The Weyl group $W$ is the subgroup of $\GL(\h^*)$ which is generated by the
reflections $s_i$ for $i\in Q_0$ such that 
\[
 s_i(\alp)= \alp - \alp(h_i)\cdot \alp_i \text{ for } \alp\in\h^*.
\]
This is a finite reflection group.

Let now $G$ be a complex simply connected simple algebraic
group with $\Lie(G)=\g$. It has maximal unipotent subgroups  $N$ 
resp. $N_-$ with $\Lie(N)=\n$ resp. $\Lie(N_-)=\n_-$, and the maximal
torus $H$ has $\Lie(H)=\h$. 
Moreover, we have standard embeddings $\vph_i\df \SL_2\ra G$ 
such that 
\[
\exp(t f_i)=\vph_i(\bsm 1&0\\t&1\esm) \text{ and } 
\exp(t e_i)=\vph_i(\bsm 1&t\\0&1\esm) \text{ for } i\in Q_0.
\]
Moreover set
\[
\eta_i(t):= \vph_i(\bsm t&0\phantom{^{-1}}\\0 &t^{-1}\esm) \in H 
\text{ for } i\in Q_0\text{ and }
t\in\CC^*.
\]
Next, recall that $N_G(H)/H$ is canonically isomorphic to the Weyl group $W$
defined above. In fact, it is possible to choose representatives 
$\bar{w}\in N_G(H)$ for the elements $w\in W$ such that
\begin{align*}
\bar{s_i}&=\exp(f_i)\exp(-e_i)\exp(f_i),\\
\overline{uv} &=\bar{u} \bar{w} \text{ if } l(uv)=l(u)+l(v).
\end{align*}
We identify the weight lattice $P=\oplus_{i\in {Q_0}}\Zet\vpi_i$
with the group of multiplicative characters of $H$  in such a way that
$\eta_i(t)^{\vpi_j}=t^{\delta_{i,j}}$ for  $i,j\in Q_0$. If we write
$?^\lam$ for the character of $H$ corresponding to the weight $\lam$ it
follows that
\[
h^{w(\lam)}=(\bar{w}^{-1}h\bar{w})^\lam\text{ for } h\in H,\lam\in P, w\in W.
\]

\subsection{Cluster algebras.} \label{LieG}
The coordinate ring of the affine base space $\CC[N_-\bcsl G]$ consists of
the functions $f\in\CC[G]$ which are invariant under $N_-$, 
i.e.~$f(g)=f(ng)$ for
all $g\in G$ and $n\in N_-$. Now $\CC[N_-\bcsl G]$ is naturally a $G$-module
via $gf(x)=f(xg)$ for $g,x\in G$. It is well-known that each irreducible
highest weight $G$-module $L(\lam)$ can be realized as a direct summand of
$\CC[N_-\bcsl G]$ by taking
\[
L(\lam) = \{f\in\CC[N_-\bcsl G]\mid f(hg)=h^\lam f(g) 
\text{ for } h\in H, g\in G\}.
\]
For each $L(\lam)$ we choose a highest weight vector $u_\lam$ which we
normalize by the condition $u_\lam(1_G)=1$. 
Following~\cite{BerZel97} we define for each fundamental weight $\vpi_i$ 
generalized minors
\[
\Del_{\vpi_i,w(\vpi_i)}:= \bar{w}\cdot u_{\vpi_i}\in L(\vpi_i)
\]
for any $w\in W$. In~\cite{BeFoZe05} it is shown in particular that the
coordinate ring of the double Bruhat cell $G^{e,w_0}=B\cap (B_-\bar{w}_0B_-)$ 
has  the structure of an (upper) cluster algebra. 
Here, $B$ and $B_-$ are opposite Borel subgroups of $G$ with $B\supset N$ and
$B_-\supset N_-$.

For a reduced expression 
$\ii=(i_1,i_2,\ldots, i_r)$ for $w_0$ which is adapted to $Q$ and
$k\in [-n,-1]\cup [1,r]$ we set
$v_{>k} := s_{i_r} s_{i_{r-1}}\cdots s_{i_{k+1}}$ if $k\geq 1$ and
$v_{>k} := w_0$ if $k\leq -1$. Then, following~\cite{BeFoZe05} set
\[
\Del(k,\ii):= \Del_{\vpi_{i_k},v_{>k}(\vpi_{i_k})}
 \]
where we take $i_k=-k$ for $k\in [-n,-1]$. 
The $\Del(k,\ii)$ for $k\in [-n,-1]\cup [1,r]$ form an initial cluster for 
$\CC[G^{e,w_0}]$. 
There is also an easy algorithm to calculate from $\ii$ the corresponding 
exchange matrix. Now set
\[
e(\ii):= \{i\in [1,r]\mid x(i)\in\Obj(\AuC{Q})\text{ non-projective}\}.
\]
There is a closely related upper cluster algebra structure
on the coordinate ring $\CC[N]$, whose initial cluster is given
by the restrictions to $N$ of the functions $\Del(k,\ii)$ with
$k\in [-n,-1]\cup e(\ii)$ (see~\ref{ssec:clustN}).
The quiver associated to the corresponding exchange matrix
is obtained from $\cA{Q}$ by removing the 
arrows between injective vertices.
See Section~\ref{sec:clustruc} for a more detailed discussion.

\subsection{Semicanonical basis.} \label{intro:scb}
In~\cite{Lusz00} Lusztig introduces the 
semicanonical basis
of the enveloping algebra $U(\n)$. Its elements are labelled naturally by
the irreducible components of the preprojective varieties
$\Lam_\dv$ for $\dv\in\ENn^{Q_0}$. In~\cite{GeLeSc04a} we started the
study of the dual semicanonical basis which can be regarded as a 
basis of $\CC[N]$.
In particular, we found that to (the isoclass of) a rigid 
$\Lam$-module $M$ there corresponds naturally an element $\rho_M$ of
this basis. If we set
\[
\theta\df [1,r]\ra [-n,-1]\cup e(\ii), j \mapsto\begin{cases}
-i &\text{ if } \RA_Q(x(j))\cong D\kQ(-,i),\\
\phantom{-}k &\text{ if } \tau^{-1}x(j)=x(k),
\end{cases}
\]
we can state our second main result precisely:
\begin{theo} \label{th1}
For $j\in [1,r]$ we have 
$\Del(j,\ii)':=\Del(\theta(j),\ii)=\rho_{\sM(j,\ii)}$. 
\end{theo}
The proof of this result is done after some preparation in 
Section~\ref{sec:cv-scb}.

\subsection{Dualities.} Denote by $?^*$ the involution on $\bQ$ with
$a\mapsto a^*$ and $a^*\mapsto a$ for all $a\in Q_1$. This induces an
anti-automorphism of $\Lam$ which we also denote by $?^*$. Thus we have
a self-duality
\[
S\df \lmd{\Lam}\longleftrightarrow \lmd{\Lam}\text{ with } SM(?)= DM(?^*).
\]
Let us define for a reduced expression $\ii=(i_1,\ldots, i_r)$ for
$w_0$ which is adapted to $Q$
\[
\sP(j,\ii):= F_\lam\Jx \AuC{Q}(x(j),-)\quad\text{ and } 
\sP_Q :=\bigoplus_{j=1}^r \sP(j,\ii).
\]
Then it is easy to see that  $\sP_Q\cong S\sM_{Q\op}$.
Thus, $\Ext^1_\Lam(\sP_Q,\sP_Q)=0$ and 
$\End_\Lam(\sP_Q)\cong\End_\Lam(\sM_{Q\op})\op$. Now, $\cA{Q\op}$ is
the Gabriel quiver of $\End_\Lam(\sM_{Q\op})\op$, see Theorem~\ref{thm:start},
and it is not hard to see that $\cA{Q\op}$ may be identified with
$\cA{Q}\op$ (recall that the same happens for Auslander-Reiten quivers:
$\AR{Q\op}\cong \AR{Q}\op$). So $\End_\Lam(\sM_{Q\op})\cong\End_\Lam(\sP_Q)\op$ and $\End_\Lam(\sM_Q)\op$ have
the same Gabriel quiver $\cA{Q}$, but they are  usually not isomorphic.

On the other hand, $S$ induces the canonical anti-automorphism $?^*$ on 
Lusztig's  algebra of constructible functions $\M\cong U(\n)$ which
commutes with the comultiplication,
see~\cite[Section 3.4]{Lusz00}. This yields by duality an automorphism
$\ome$ of the coordinate ring $\CC[N]$ which anti-commutes with the
comultiplication. The corresponding anti-automorphism  of $N$ leaves
the one-parameter subgroups $\CC\ra N, t\mapsto\exp(t e_i)$ invariant.

Note that
$\ii^*:= (\mu(i_r),\mu(i_{r-1}),\ldots,\mu(i_1))$ is a reduced expression for
$w_0$ which is adapted to $Q\op$, see~\ref{ssec:cov-constr} for the definition 
of $\mu$. We leave it to the reader to verify the following:
\[
\rho_{\sP(j,\ii)}= \omega(\rho_{\sM(r+1-j,\ii^*)})\quad\text{for } 
1\leq j\leq r,
\]
see also~\ref{intro:scb}.

\subsection{Triangulated structures.} In the appendix 
(Section~\ref{sec:apx}) we point out that 
the category of injective $\tLam$-modules is triangulated. This implies that
the category of injective $\Lam$-modules is also triangulated. This fact
helps us  to explain the unusual symmetries in the stable module categories
over both categories by  results of Freyd and Heller. 
Among other useful formulas we
can recover the famous ``6-periodicity'' for modules over the preprojective
algebra.

\section{The universal cover of a preprojective algebra}
\subsection{Quiver categories.}\label{ssec:quica}
 Let $Q=(Q_0,Q_1,t,h)$ be a quiver with vertices $Q_0$,
arrows $Q_1$ and $t,h\df Q_1\ra Q_0$ such that we have $a\df t(a)\ra h(a)$
for each arrow $a\in Q_1$. A {\em path} in $Q$ is a sequence of arrows
$a_n a_{n-1}\cdots a_1$ such that $t(a_{i+1})=h(a_i)$ for $i=1,2,\ldots, n-1$.

On $\Zet^{Q_0}$ we have the Ringel bilinear form 
\[
\ebrace{\dv,\dw}=\sum_{i\in Q_0}\dv(i)\dw(i)-\sum_{a\in Q_1}\dv(t(a))\dw(h(a)).
\]
 Let $\ka$ be a field. Since we need to 
consider infinite coverings of preprojective algebras we have to 
consider $\ka$-categories rather than $\ka$-algebras. 

If $Q$ is a quiver we denote by $\ka[Q]=\kQ$ the $\ka$-category which has 
$Q_0$ as objects and with the morphism space $\kQ(p,q)$ having the paths
from $p$ to $q$ as a basis. The composition is naturally induced from
the concatenation of paths.

\subsection{Conventions.} \label{ssec:a-conv}
If $\cD$ is a $\ka$-category we denote by $\lmd{\cD}$ the category
of  finitely presented (covariant) $\ka$-functors 
$\cD\ra\lmd{\ka}$. These functors are also called left modules,
see for example~\cite[Section 2.2]{garo92} for more details.

Let  $G$ be a group of $\ka$-automorphisms of $\cD$ which acts from the
left on $\cD$. Then $G$ acts naturally from the right on $\lmd{\cD}$. 
If $g\in G$ and $M$ is a left $\cD$-module, then we write
$M^g(-):= M(g^{-1}-)$ for the twisted module. For example, if $x\in\cD$
we get for the projective module $\cD(x,-)$ an isomorphism
$\cD^g(x,-)\cong\cD(gx,-)$.

\subsection{Construction of $\tLam$.}
\label{ssec:cov-constr}
Let $Q=(Q_0,Q_1,t,h)$ be a (connected) Dynkin quiver, so the underlying 
graph  $\abs{Q}$ is one of the following:
\[
\def\objectstyle{\scriptstyle}
\xymatrix@-1.2pc{
{\textstyle{\Aa_n}}\df & 1\ar@{-}[r]& 2\ar@{-}[r]& 3\ar@{-}[r] &
\quad\cdots\quad\ar@{-}[r] & (n-2)\ar@{-}[r] & (n-1)\ar@{-}[r] &n\\
&&&&&&(n-1)\\
{\textstyle{\Dd_n}}\df & 1\ar@{-}[r]& 2\ar@{-}[r]& 3\ar@{-}[r] &
\quad\cdots\quad\ar@{-}[r] & (n-2)\ar@{-}[ru]\ar@{-}[rd] & &&{n\geq 4}\\
&&&&&&n\\
{\textstyle{\Ee_n}}\df & 1\ar@{-}[r]& 2\ar@{-}[r]& 3\ar@{-}[r]\ar@{-}[d] &
\quad\cdots\quad\ar@{-}[r] & (n-2)\ar@{-}[r] & (n-1)&&n=6,7,8\\
&&&n}
\]
Define a (possibly trivial) involution $\mu$ on the vertices of $Q$ by
\[
\mu(q)=
\begin{cases}
n+1-q  &\text{ in case } \Aa_n,\\
2n-1-q &\text{ in case } \Dd_n \text{ and } (n \text{ odd and } q\geq n-1),\\
6-q    &\text{ in case } \Ee_6 \text{ and } q\leq 5,\\
q      &\text{ otherwise.} 
\end{cases}
\]
Following~\cite{gabr80} let
$\ZQ$ be the quiver with vertices $\ZQ_0=\Zet\times Q_0$ and
arrows $\ZQ_1=\Zet\times\{Q_1\cup Q_1^*\}$ where
\begin{alignat*}{2} 
\tilde{t}(i,a)  &=(i+1,t(a)),&\qquad
\tilde{t}(i,a^*)&=(i,h(a)),\\ 
\tilde{h}(i,\alp)  &=(i,h(a)),&\qquad
\tilde{h}(i,a^*)&=(i,t(a)).
\end{alignat*}
We  think of $Q\op$ as a subquiver of $\ZQ$ 
via the embedding $q\mapsto (0,q)$. 
The quiver $\ZQ$ admits a translation automorphism $\tau$ induced by
$\tau(i,q)=(i+1,q)$. Moreover we have a ``Nakayama permutation''
$\hnu$ of $\ZQ$. In order to define it we need  
the following auxiliary function: For any vertex $q$ denote by
$l(q)$ the number of arrows pointing towards $1$ on the unique walk from $1$ 
to $q$. Now, $\hnu$  is defined on the vertices by
\[
\hnu(p,q)=\begin{cases}
(p+q-1+l(\mu(q))-l(q),\mu(q))&\text{ in case }\Aa_n,\\
(p+n-2+l(\mu(q))-l(q),\mu(q))&\text{ in case }\Dd_n,\\
(p+q+2+l(\mu(q))-l(q),\mu(q))&\text{ in case }\Ee_6 \text{ and } q\leq 5,\\
(p+5,\mu(q))  &\text{ in case }\Ee_6 \text{ and } q=6,\\
(p+8,\mu(q))  &\text{ in case }\Ee_7,\\
(p+14,\mu(q)) &\text{ in case }\Ee_8.
\end{cases}        
\]
So, $\hnu$ is a ``translation reflection'' stabilizing the ``middle line'' in
case $\Aa_n$ and a translation in the cases $\Dd_{2k}$, $\Ee_7$ and $\Ee_8$.
At the beginning of~\ref{expl-GaD5}
we show as an example part of $\ZQ$ for a quiver
of type 
$\Dd_5$, together with the action of $\hnu$.

We consider the {\em mesh ideal} $I$ in the category $\ka[\ZQ]$ 
which is generated by the elements
\begin{equation} \label{eqn:mesh-rel}
\sum_{\substack{a\in Q_1\\ t(a)=q}} (i,a^*)\cdot (i,a) -
\sum_{\substack{a\in Q_1\\ h(a)=q}} (i,a)\cdot (i+1,a^*)
\quad\text{ for all } (i,q)\in\ZQ_0.
\end{equation}
Note that $\tLam:=\ka[\ZQ]/I$ is independent of the orientation
of $Q$ (up to relabelling of the vertices), and that $I$ is invariant
under the quiver automorphisms $\tau$ and $\hnu$. Thus we will use the same
names for the induced automorphisms of $\tLam$. Note further that 
the commuting automorphisms $\tau$ and $\hnu$ of $\tLam$ are determined up to 
isomorphism by their effect on objects.

The action of the group $\Zet$ via $\tau$ 
on $\tLam$ induces self-equivalences $?^{(i)}$ of $\lmd{\tLam}$ with 
$M^{(i)}(-):= M(\tau^{-i} -)$ for $i\in\Zet$. Moreover we obtain
the covering functor
$F\df\tLam\ra\Lam$ which sends $(i,q)$ to $q$. Associated to $F$ we have the
push-down $F_\lam\df\lmd{\tLam}\ra\lmd{\Lam}$ with 
\[
(F_\lam M)(q)=\oplus_{i\in\Zet} M(i,q)
\]
and the obvious effect on morphisms.

\subsection{Auslander category.} \label{ssec:ausca}
We define a function 
$N\df Q_0\ra\ENn$ by the property
\[
\tau^{N(q)}(0,q)=\hnu(0,\mu(q)).
\]
This is well-defined by the construction of $\hnu$ since 
$\mu$ is an involution on $Q_0$, see~\ref{ssec:cov-constr}.
The function $N$  depends on the orientation of $Q$ in case $\Aa_n$, 
$\Dd_{2k+1}$ and $\Ee_6$. In any case we have $N(q)+N(\mu(q))=h(Q)-2$, where
$h(Q)$ denotes the Coxeter number of $\abs{Q}$.

Define now $\AuC{Q}$ as the full subcategory of $\tLam$ which has the objects
of the form
$(i,q)=\tau^{i-N(q)}\hnu(0,\mu(q))$ with $q\in Q_0$ and $0\leq i\leq N(q)$.
In other words, we take the objects which lie between the two copies of
$Q\op$ in $\ZQ$ which are obtained via $q\mapsto (0,q)$ 
resp.~ $q\mapsto\hnu(0,q)$.
This is the {\em Auslander category} of $\kQ$. Note that $\AuC{Q}$ depends
on the orientation of $Q$.

By construction we have a full embedding $\iot=\iot_Q\!\df\kQ\op\ra\AuC{Q}$ 
induced by $q\mapsto (0,q)$. 
Moreover, $\AuC{Q}$ is canonically equivalent to the category
of indecomposable $\kQ$-modules via the functor
\begin{equation} \label{eqn:ausiso}
\RA_Q\df\AuC{Q}\ra\lind{\kQ},\ x\mapsto \AuC{Q}(\hnu\iot,x),
\end{equation}
where $\AuC{Q}(\hnu\iot,x)=\AuC{Q}(-,x)\circ\hnu\iot\df\kQ\op\ra\lmd{\ka}$ 
is a  contravariant functor which we have to interpret as a left $\kQ$-module.
For example, $\RA_Q(0,q)\cong D\kQ(-,q)$ and 
$\RA_Q(\hnu(0,q))\cong\kQ(q,-)$.

Thus the Gabriel quiver of $\AuC{Q}$ (which is the full subquiver of
$\ZQ$ with the vertices from $\AuC{Q}$) is the Auslander-Reiten quiver $A_Q$
of $\kQ$.

Similarly, since  we consider left modules, one obtains an equivalence
\[
\linj{\AuC{Q}}\ra\lmd{\kQ\op}\cong(\lmd{\kQ})\op, 
I\mapsto I\circ\hnu\iot.
\]
Here, $\linj{\AuC{Q}}$ denotes the category of injective left 
$\AuC{Q}$-modules.

Note that  the $q\in Q_0$ parametrize the  indecomposable
projective-injective $\AuC{Q}$-modules, namely we have
\[
D\AuC{Q}(-,(0,q))\cong\AuC{Q}(\hnu(0,q),-).
\]

Recall, that as an Auslander category $\AuC{Q}$ has dominant dimension at 
least~$2$~\cite[VI.5]{AuReSm97}.
This means by duality that each (indecomposable) injective $\AuC{Q}$-module 
$D\AuC{Q}(-,x)$ has a projective presentation
\[
P_{1,x}\ra P_{0,x}\ra D\AuC{Q}(-,x)\ra 0
\]
with $P_{0,x}$ and $P_{1,x}$ also injective. 

\section{Start modules} \label{sec:sm}
\subsection{Adjoint functors.} \label{ssec:adj}
Let $J\df\AuC{Q}\ra\tLam$ be the full embedding of 
locally bounded categories. Since $\AuC{Q}$ is convex in $\tLam$ we get an 
exact functor ``extension by $0$''
\[
\Jx\df  \lmd{\AuC{Q}}\ra\lmd{\tLam}\text{ with }
(\Jx M)(x)=
\begin{cases} M(x) &\text{ if } x\in\AuC{Q},\\ 0 &\text{ else.}\end{cases}
\]
For a $\AuC{Q}$-module $N$ and $x\in\AuC{Q}$ we have natural isomorphisms
\[
\Hom_{\AuC{Q}}(N,D\AuC{Q}(-,x))\cong DN(x)\cong\Hom_{\tLam}(\Jx N,D\tLam(-,x)).
\]
We conclude that $\Jx$ has a right  adjoint $\Jr$ which is defined by
being left exact and
\[
\Jr D\tLam(-,x)=
\begin{cases} D\AuC{Q}(-,x)&\text{ if } x\in\Obj(\AuC{Q}),\\ 
                         0&\text{ else. }
\end{cases}
\]
The adjunction morphism $\iot_M\colon  \Jx \Jr M\ra M$ is
injective, its co-kernel is co-generated by 
$\bigoplus_{y\in\Obj(\tLam)\setminus\Obj(\AuC{Q})} M(y)$.

\begin{lem} \label{lem:aus}
Let $x,y\in\Obj(\AuC{Q})$ and $i\in\Zet$. Then
\[
\Hom_{\tLam}(\Jx D\AuC{Q}(-,x),\Jx D\AuC{Q}^{(i)}(-,y))\cong
\begin{cases}
\AuC{Q}(\tau^i y,x) &\text{ if }i\geq 0\text{ and }\tau^i y\in\Obj(\AuC{Q}),\\
0                &\text{ else.}
\end{cases}
\]
In the first case we have more generally 
$\Jr(\Jx D\AuC{Q}^{(i)}(-,y))\cong D\AuC{Q}(-,\tau^i y)$.
\end{lem}

Note, that $\Jx D\Gam(-,y)$ is the injective $\Gam_Q$-module with socle 
concentrated in $y$, but seen as $\tLam$-module, thus we may apply to it
the translation functor $?^{(i)}$, see~\ref{ssec:cov-constr}.

\begin{pf}
The $\tLam$-module $M={\Jx D\AuC{Q}^{(i)}(-,y)}$ has simple socle concentrated 
in the one-dimensional space $M(\tau^i y)$. If $\tau^i y\not\in\Obj(\AuC{Q})$ 
then this does not belong to the support of $\Jx D\AuC{Q}(-,x)$.
On the other hand, if $i<0$, it is sufficient to show that there are no
maps from an induced projective-injective  module to $M$ since $\AuC{Q}$ has
dominant dimension $\geq 2$, see~\ref{ssec:ausca}. Now,
\[
0=M(\hnu(0,q))=\Hom_{\tLam}(\tLam(\hnu(0,q),-),M)=
\Hom_{\tLam}(\Jx \AuC{Q}(\hnu(0,q),-),M)
\]
with the first identity holding for $i<0$.

Thus, let $i\geq 0$ and $\tau^i y\in\AuC{Q}$. In this case we have in
$\lmd{\tLam}$ an injective presentation
\[
0 \ra M \ra D\tLam(-,\tau^i x)
   \ra \oplus_j D\tLam(-,y_j)^{m(j)}
\]
for certain $y_j\in\Obj(\tLam)\setminus\Obj(\AuC{Q})$. 
Our claim follows now from the construction of $\Jr$.
\end{pf}

\subsection{A graded category.}
\label{ssec:start-grad}
We construct the  $\ENn$-graded category $\cGam_Q$. It has the same objects as
$\AuC{Q}$, but the homogenous components are 
$\cGam_{Q,i}(x,y):=\AuC{Q}(\tau^i x,y)$ if $\tau^ix\in\Obj(\AuC{Q})$ and
$\cGam_{Q,i}(x,-)=0$ otherwise. 
The natural composition is given by
\[
\cGam_{Q,j}(y,z)\otimes\cGam_{Q,i}(x,y)\ra\cGam_{Q,i+j}(x,z),\quad
(\psi\otimes \phi)\mapsto 
\begin{cases}\psi\circ(\tau^j\phi) &\text{ if } \tau^{i+j}x\in\Obj(\AuC{Q}),\\
             0,                    &\text{ else.}\end{cases}
\]
By construction, each morphism in $\cGam_{Q,\geq 1}(x,-)$ factors through
$\id_{\tau x}\in\cGam_{Q,1}(x,\tau x)$ if $\tau x\in\Obj(\AuC{Q})$, 
otherwise $\cGam_{Q,\geq 1}(x,-)=0$. 
Moreover we have
\[
\id_{\tau^n x}\circ \cdots \circ
\id_{\tau^{2} x}\circ \id_{\tau x}=\id_{\tau^n x}\in\cGam_{Q,n}(x,\tau^n x)
\text{ if }\tau^{n}x\in\Obj(\AuC{Q}),
\]
where we consider on the left hand side 
$\id_{\tau^i x}\in \cGam_{Q,1}(\tau^{i-1}x,\tau^i x)$.
Now, $\cGam_Q$ can be described easily by a graded quiver $\cA{Q}$. It has the
same vertices and degree 0 arrows as the Auslander-Reiten quiver $\AR{Q}$ of
$\kQ$ (i.e. the quiver of $\AuC{Q}$) moreover
there is a degree 1 arrow $t_x\colon x\ra \tau x$ for each
$x\in\Obj(\AuC{Q})$ with $\tau x\in\Obj(\AuC{Q})$. The degree $0$ 
relations are the mesh relations for $\AuC{Q}$. 
Moreover, each degree $0$ arrow 
$a\df x\ra y$ with  $y$ not projective gives rise
to a degree $1$ relation
\[
t_y a-(\tau a)t_x.
\]
This has to be interpreted as a zero-relation if $x$  is projective.
A nice way to remember these relations is the following:
For each arrow between to vertices which are not both  injective
there is a (generic) homogeneous length 2 relation in the opposite direction, 
see also~\ref{expl-ga*}.

\subsection{Dynkin quivers.}
\label{ssec:start-prep}
Let us collect some basic facts about the representation theory of 
a Dynkin quiver $Q$. Define $\di_q:=\dimv D\kQ(-,q)$
and $\dpr_q:=\dimv\kQ(q,-)$ for $q\in Q_0$, the dimension vectors of the
indecomposable injective and projective $\kQ$-modules, respectively. 
We have then
for $0\leq i\leq N(q)$
\begin{equation}\label{eq:cox-tau}
\Phi^i\di_q=\dimv(\tau_Q^i D\kQ(-,q))=\dimv(\tau_Q^{i-N(q)}\kQ(q,-))=
\Phi^{i-N(q)}\dpr_{\mu(q)},
\end{equation}
where $\Phi$ denotes the Coxeter transformation and $\tau_Q$ the 
Auslander-Reiten translate in $\lmd{\kQ}$.
Next, if $\ebrace{-,-}$ denotes the Ringel bilinear form of
$\kQ$ we have
\begin{itemize}
\item
$\ebrace{\dimv M,\dimv N}=\dim\Hom_Q(M,N)-\dim\Ext^1_Q(M,N)$,
\item
$\ebrace{\dv,\dw}=\ebrace{\Phi\dv,\Phi\dw}=-\ebrace{\dw,\Phi\dv}$,
\item
$\ebrace{\dv,\di_q}=\dv(q)$ thus $\ebrace{\dv,\sum_{q\in Q_0}\di_q}=\abs{\dv}$,
where $\abs{\dv}:=\sum_{q\in Q_0}\dv(q)$,
\end{itemize}
see for example~\cite[2.4]{ring84}.

\begin{lemu} \label{lemu:dim}
Let $\tau^i(0,p)=(i,p)$ and $\tau^{-j}\hnu(0,q)$ belong to $\Obj(\AuC{Q})$, 
then 
\[
\dim\AuC{Q}(\tau^{-j}\hnu(0,q),\tau^i(0,p))=
\begin{cases}
(\Phi^{i+j}\di_p)(q)  &\text{ if } i+j\leq N(p),\\
(\Phi^{-i-j}\dpr_q)(p)&\text{ if } i+j\leq N(\mu(q)),\\
0                     &\text{ if } i+j>\min\{N(p),N(\mu(q))\}.
\end{cases}
\] 
\end{lemu}
Note that the three cases are not exclusive, however they cover
obviously all possibilities for $(i,j)\in [0,N(p)]\times [0,N(\mu(q))]$.

\begin{pf}
In the first case we have 
\begin{equation} \label{eq:lemu-dim}
\dim\AuC{Q}(\tau^{-j}\hnu(0,q),\tau^i(0,p))=
\dim\AuC{Q}(\hnu(0,q),\tau^{i+j}(0,q)).
\end{equation}
On the other hand, the equivalence $\RA_Q$ from $\AuC{Q}$ to the category of 
indecomposable $\kQ$-modules commutes with translations
\[
\RA_Q(\tau^i(0,q))\cong \tau_Q^i D\kQ(-,q)\quad\text{for } 0\leq i\leq N(q).
\]
Thus~\eqref{eq:lemu-dim} is equal to
\[
\dim\Hom_Q(\kQ(q,-),\tau_Q^{i+j}D\kQ(-,p))
=\dim\Hom_Q(\tau_Q^{-j}\kQ(q,-),\tau_Q^iD\kQ(-,p)).
\]
Our claim follows now from~\eqref{eq:cox-tau} since  for a 
finite-dimensional $\kQ$-module $M$ we have that
\[
\Hom_Q(\kQ(q,-),M)\cong M(q)\cong D\Hom_Q(M,D\kQ(-,q)).
\]
The second case is treated similarly.
Finally we have
\[
\AuC{Q}(\tau^{-j}\hnu(0,q),\tau^i(0,p))\cong 
\tLam(\tau^{-j}\hnu(0,q),\tau^i(0,p))
\cong\tLam(\hnu(0,q),\tau^{i+j-N(p)}\hnu(0,\mu(p))).
\]
The last term vanishes obviously for $i+j>N(p)$. A similar argument shows that
$\AuC{Q}(\tau^{-j}\hnu(0,q),\tau^i(0,p))=0$ for $i+j>N(\mu(q))$.
\end{pf}

\begin{lemu} \label{lemu:trick}
If $N(p)-N(q)> j\geq 0$ holds for some $p,q\in Q_0$, 
then $\ebrace{\Phi^j\di_p,\di_q}=0$.
\end{lemu}

\begin{pf} Since $D\kQ(-,q)$ is injective we have
\[
\ebrace{\Phi^j\di_p,\di_q}=\dim\Hom_Q(\tau_Q^j D\kQ(-,p),D\kQ(-,q))=
\dim\tLam(\tau^j(0,p),(0,q)).
\] 
On the other hand
for $N(p)-j>N(q)$ there is no path from
$\tau^j(0,p)=(N(p)-j,\mu(p))$ to $(0,q)$ in $\Zet Q$.
\end{pf}

\begin{prp} \label{prp:form}
With the notation of~\ref{ssec:startM} and~\ref{ssec:start-prep} we have
\begin{equation} \label{eq:dimv}
(\dimv \sM_Q)(\mu(p)) =\sum_{q\in Q_0}\sum_{i=0}^{N(q)}(i+1)(\Phi^i \di_q)(p)
\end{equation}
for $p\in Q_0$, and
\begin{equation} \label{eq:endim}
\dim\cE =\sum_{q\in Q_0}\sum_{i=0}^{N(q)} 
\left( \binom{N(q)+2}{2}-\binom{i+1}{2}\right)\abs{\Phi^i \di_q}.
\end{equation}
\end{prp}

\begin{pf}
For~\eqref{eq:dimv} we observe first that
\[
\dimv F_\lam \Jx D\AuC{Q}(-,(j,q))(\mu(p))=\sum_{i=j}^{N(q)}(\Phi^i\di_q)(p)
\]
for $(j,q)\in\Obj(\AuC{Q})$ by~\ref{lemu:dim} and the definition of the 
push-down $F_\lam$. Now~\eqref{eq:dimv} follows from the definition of $\sM_Q$.

For~\eqref{eq:endim} we observe first that
\[
\abs{\dimv F_\lam \Jx \AuC{Q}(\tau^{-j}\hnu(0,q),-)}=
\sum_{k=j}^{N(\mu(q))}\abs{\Phi^{-k}\dpr_q}
\]
for $\tau^{-j}\hnu(0,q)\in\Obj(\AuC{Q})$ again by~\ref{lemu:dim} and the 
definition of the push-down. Now, by construction of $\cE$ we have
\begin{align*}
\dim\cE &= \sum_{q\in Q_0}\sum_{i=0}^{N(\mu(q))}\sum_{j=0}^i
\abs{\dimv F_\lam\Jx\AuC{Q}(\tau^{-j}\hnu(0,q),-)}\\
&=\sum_{q\in Q_0}\sum_{i=0}^{N(\mu(q))}\sum_{j=0}^i\sum_{k=i}^{N(\mu(q))}
\abs{\Phi^{-k}\dpr_q}\\
&=\sum_{q\in Q_0}\sum_{i=0}^{N(q)}\sum_{j=0}^i\sum_{k=i}^{N(q)-j}
\abs{\Phi^k\di_q}
=\sum_{q\in Q_0}\sum_{i=0}^{N(q)}\sum_{j=i}^{N(q)+1} j\abs{\Phi^i\di_q}.
\end{align*}
\end{pf}

\begin{prp} \label{prp:start-eqn}
Let $Q$ be a Dynkin quiver. Then for
$\dv= \sum_{q\in Q_0}\sum_{d=0}^{N(q)} (i+1)\Phi^i\di_q$  holds 
\[
\ebrace{\dv,\dv}=\sum_{p\in Q_0}\sum_{d=0}^{N(p)}
\left(\binom{N(p)+2}{2}-\binom{d+1}{2}\right)\abs{\Phi^d\di_p}.
\]
\end{prp}
\begin{pf}
For  convenience let us write
$N(p,q):=\{0,1,\ldots, N(p)\}\times\{0,1,\ldots, N(q)\}$ and
${N(p,q,\geq)}:={\{(i,j)\in N(p,q)\mid i\geq j\}}$, similarly
${N(p,q,<)}:= {N(p,q)\setminus N(p,q,\geq)}$.
Now we have
\begin{align*}
\ebrace{\dv,\dv}= &\phantom{-}\sum_{p,q\in Q_0}\quad
\sum_{(i,j)\in N(p,q)}(i+1)(j+1)\ebrace{\Phi^i\di_p,\Phi^j\di_q}\\
= & \phantom{-}\sum_{p,q\in Q_0}\quad\!\! \sum_{(i,j)\in N(p,q,\geq)}(i+1)(j+1)
\ebrace{\Phi^i\di_p,\Phi^j\di_q}\\
&\!-\!\sum_{p,q\in Q_0}\quad\!\!\sum_{(i,j)\in N(p,q,<)} (i+1)(j+1)
\ebrace{\Phi^j\di_q,\Phi^{i+1}\di_p} \displaybreak[0]\\ 
=& \phantom{-}\sum_{p,q\in Q_0}\quad\! \sum_{(i,j)\in N(p,q,\geq)} (i+1)
\ebrace{\Phi^i\di_p,\Phi^j\di_q}\\
&\!-\!\sum_{p,q\in Q_0}\quad\;\;\sum_{i=N(q)+1}^{N(p)}(i+1)(N(q)+1)
\ebrace{\Phi^i\di_p,\Phi^{N(q)+1}\di_q}\\
\intertext{here, the second sum vanishes by Lemma~\ref{lemu:trick}, thus 
from~\ref{ssec:start-prep}}
=&\phantom{-}\sum_{p,q\in Q_0}\quad\sum_{(i,j)\in N(p,q,\geq)}(i+1)
\ebrace{\Phi^{i-j}\di_p,\di_q}\\
=&\phantom{-}\sum_{p,q\in Q_0}\quad\sum_{d=0}^{N(p)}
\left(\sum_{k=d}^{\min\{N(p),N(q)+d\}}(k+1)\right)
\ebrace{\Phi^d\di_p,\di_q}\\
=&\phantom{-}\;\sum_{p\in Q_0}\quad \sum_{d=0}^{N(p)}\sum_{k=d}^{N(p)}(k+1)
\ebrace{\Phi^d\di_p,\sum_{q\in Q_0}\di_q}\\
&-\sum_{p,q\in Q_0}\quad\sum_{d=0}^{N(p)}\;\;\sum_{k=N(p)+d}^{N(p)}(k+1)
\ebrace{\Phi^d\di_p,\di_q}.
\end{align*}
Here again, the second sum vanishes by Lemma~\ref{lemu:trick}. 
Finally, $\ebrace{\Phi^d\di_p,\sum_{q\in Q_0}\di_q}=\abs{\Phi^d\di_p}$ as observed 
at the beginning of~\ref{ssec:start-prep}, so our claim follows.
\end{pf}

\subsection{Proof of Theorem~\ref{thm:start}.} \label{ssec:start-pf}
The first claim follows directly from the construction of $\cGam_Q$,
Lemma~\ref{lem:aus} and the fact that
\[
\Hom_\Lam(F_\lam N, F_\lam N)\cong\oplus_{i\in\Zet}\Hom_{\tLam}(N,N^{(i)}).
\]
For the second claim we  have to show by~\cite[Lemma 1]{cb00} that
$\ebrace{\dimv \sM_Q,\dimv \sM_Q}=\dim\mathcal{E}$.
This   follows from the dimension formulas~\eqref{eq:dimv} 
and~\eqref{eq:endim} of Proposition~\ref{prp:form} together with
Proposition~\ref{prp:start-eqn}.

\section{The cluster algebras $\CC[G^{e,w_0}]$ and $\CC[N]$}
\label{sec:clustruc}
In the next two sections we will use the setup 
from~\ref{intro:Liea} and~\ref{LieG}.
We will need the following  result, see for 
instance~\cite[Section 4.4.3]{Joseph95}

\begin{lem} \label{lem:well-known}
Let $\ii=(i_1,i_2,\ldots, i_m)$ be a reduced expression for
some element $w^{-1}\in W$, and $L(\lam)$ an irreducible representation of
highest weight $\lam$ for $G$. If $u_\lam$ is a highest weight vector
for $L(\lam)$ then we have
\[
\bar{w} u_\lam=f_{i_m}^{(b_m)}\cdots f_{i_2}^{(b_2)}f_{i_1}^{(b_1)} (u_\lam)
\text{ and } f_{i_m} (\bar{w} u_\lam)=0.
\]
Here $b_1:= \lam(h_{i_1})$,  
$b_k:=(s_{i_{k-1}}\cdots s_{i_2} s_{i_1}(\lam))(h_{i_k})$ for $2\leq k\leq m$
and 
$f_{i_k}^{(b_k)}:=\frac{1}{b_k!}f_{i_k}^{b_k}\in U(\g)$.
\end{lem}

For $v\in L(\lam)$, we write 
$f_j^{\max}(v) := (1/m!)f_j^m(v)$ where $m=\max\{p \mid f_j^p(v) \not = 0\}$.
With this notation we could restate the equality of the lemma as
\[
\bar{w} u_\lam=f_{i_m}^{\max}\cdots f_{i_2}^{\max}f_{i_1}^{\max} (u_\lam).
\]

\subsection{} \label{1;2}
The group $G$ has Bruhat decompositions with respect to $B$ and
$B_{-}$, namely
\[
G=\bigcup_{u\in W} B\bar{u}B = \bigcup_{v\in W} B_{-}\bar{v}B_{-}.
\]
The intersection of two cells $G^{u,v} = (B\bar{u}B) \cap (B_{-}\bar{v}B_{-})$
is called a {\em double Bruhat cell}. 

In particular taking $u=e$, the
unit in $W$, and $v=w_0$ we obtain
$G^{e,w_0} = B \cap (B_{-}\bar{w}_0B_{-})$, the intersection of $B$ with
the big cell relative to $B_{-}$.
By~\cite[Proposition 2.8]{BeFoZe05}, $G^{e,w_0}$
consists of all elements $x$ of $B$ such that 
$\Delta_{\vpi_i,w_0(\vpi_i)}(x)\not = 0$ for every $i$.
This is a Zariski open subset of $B$, hence an algebraic variety 
of dimension $n+r$, where $r=|\pRt|$ is the number of positive roots
associated to the Dynkin type of $G$.
Moreover, we see that the algebra of regular functions $\CC[G^{e,w_0}]$ 
is obtained from $\CC[B]$ by adjoining formal inverses to the 
functions $\Delta_{\vpi_i,w_0(\vpi_i)}$. 

On the other hand, $N$ can be described as the subvariety of $B$
given by the equations 
\[
\Delta_{\vpi_i,\vpi_i}(x) = 1, \qquad (1\le i\le n).
\]
Hence the algebra $\CC[N]$ is the quotient of $\CC[B]$ by the ideal generated 
by the elements $(\Del_{\vpi_i,\vpi_i}-1)_{i=1,2,\ldots,n}$.

\subsection{}
In \cite{FomZel99}, Fomin and Zelevinsky have introduced a transcendence
basis $F(\ii)$ of the field of rational functions $\CC(G^{e,w_0})$
consisting of certain generalized minors.
In \cite{BeFoZe05}
Berenstein, Fomin and Zelevinsky have shown that
each $F(\ii)$ can be taken as the initial cluster for 
a natural upper cluster algebra structure on the ring
$\CC[G^{e,w_0}]$.
We are now going to recall their construction. 

\subsubsection{}
We add $n$ additional letters 
$i_{-n},\ldots,i_{-1}$ at the beginning of $\ii$, where
$i_{-j}=-j$, and obtain an $(r+n)$-tuple
\[ 
(i_{-n},\ldots,i_{-1},i_1,\ldots, i_r)=(-n,\ldots,-1,i_1,\ldots, i_r). 
\]
For $k \in [-n,-1] \cup [1,r]$ let 
\[
k^+ =  
\begin{cases}
r+1 & \text{if $|i_l| \not= |i_k|$ for all $l > k$},\\
\min\{l \mid l > k \text{ and } |i_l| = |i_k| \} & \text{otherwise}.
\end{cases}
\]
Then $k$ is called $\ii$-{\it exchangeable} if $k$ and $k^+$ are
in $[1,r]$.
Let $e(\ii)\subset [1,r]$ be the set of $\ii$-exchangeable 
elements.
One easily checks that $e(\ii)$ contains $r-n$ elements.
More precisely, the set of indices $i_k$ for $k\in [1,r] - e(\ii)$
is exactly $[1,n]$.

\subsubsection{}
Next, one defines a quiver $\tAi$ with set of vertices
$[-n,-1] \cup [1,r]$.
Assume that $k$ and $l$ are vertices such that the following 
hold:
\begin{itemize}

\item $k < l$;

\item $\{ k,l \} \cap e(\ii) \not= \emptyset$.

\end{itemize}
%
There is an arrow $k \to l$ in $\tAi$ if and only if
$k^+ = l$, 
and there is an arrow $l \to k$ if and only if
$l < k^+ < l^+$ and $a_{|i_k|,|i_l|} = -1$.
Here, $(a_{ij})_{1\le i,j \le n}$ denotes the Cartan matrix
of the root system of $G$.
By definition these are all the arrows of $\tAi$.

\begin{remu} \label{ex1}
If $\ii$ is a reduced expression for $w_0$ which is adapted to a Dynkin quiver $Q$,
then it is easy to obtain $\tAi$ from the Auslander-Reiten
quiver $A_Q$. See the examples in~\ref{expl:ada-ord} and~\ref{expl:Gam-i}. 
\end{remu}

\subsubsection{}
Now define an $(r+n) \times (r-n)$-matrix 
\[
\widetilde{B}(\ii) = (b_{kl})
\]
as follows.
The columns of $\widetilde{B}(\ii)$ are indexed by the elements
in $e(\ii)$, and the rows by $[-n,-1] \cup [1,r]$.
Set
\[
b_{kl} =
\begin{cases}
\phantom{-}1&\text{if there is an arrow $k\to l$ in $\tAi$},\\
-1 & \text{if there is an arrow $l \to k$ in  $\tAi$},\\
\phantom{-}0 & \text{otherwise}.
\end{cases}
\]

\subsubsection{} \label{134}
For $k\in [-n,-1] \cup [1,r]$ one defines a generalized minor
$\Delta(k,\ii)$ as follows.
For $k\in [1,r]$ set 
$
v_{>k} = s_{i_r} s_{i_{r-1}} \cdots s_{i_{k+1}}
$
and for $k\in [-n,-1]$ put $v_{>k} = w_0$.
Then define
\[
\Delta(k,\ii) = \Delta_{\vpi_{|i_k|},\, v_{>k}(\vpi_{|i_k|})}.
\]
Since $s_j(\vpi_i)=\vpi_i$ for $j \not = i$, it is easy to see
that if $k \in [1,r]$ is not exchangeable then 
$\Delta(k,\ii) = \Delta_{\vpi_{i_k},\, \vpi_{i_k}}$.
On the other hand for $-i\in [-n,-1]$ we have 
$\Delta(-i,\ii) = \Delta_{\vpi_i,\, w_0(\vpi_i)}$.

It is known \cite[Theorem 1.12]{FomZel99} that this collection of
$n+r$ minors is a transcendence basis of the field $\CC(G^{e,w_0})$,
for any reduced expression $\ii$ of $w_0$.
By \ref{1;2}, we see that if we remove from this collection
the $n$ minors $\Delta_{\vpi_i,\,\vpi_i}$ we obtain a
transcendence basis of the field $\CC(N)$.

\subsubsection{} \label{1;3;4}
Let $\mathcal{F}$ be the field of rational functions over $\CC$
in $n+r$ independent variables 
$\tilde{\mathbf{x}}=(x_{-n},\ldots,x_{-1},x_1,\ldots,x_r)$.
Let $\overline{\mathcal{A}}(\ii)_\CC$ denote the upper cluster algebra
associated to the seed $(\tilde{\mathbf{x}},\widetilde{B}(\ii))$,
a subalgebra of $\mathcal{F}$ (see \cite[Definition 1.6]{BeFoZe05}).
Here the non-exchangeable indices in $[-n,-1] \cup [1,r]$ label
the generators of the coefficient group (see \cite[\S 2.2]{BeFoZe05}).

Berenstein, Fomin and Zelevinsky then show that the isomorphism of fields 
$\varphi_\ii$ from $\mathcal{F}$ to $\CC(G^{e,w_0})$ defined by
\[
\varphi_\ii(x_k)=\Delta(k,\ii), \qquad (k\in [-n,-1]\cup[1,r]),
\] 
restricts to an algebra isomorphism 
$\overline{\mathcal{A}}(\ii)_\CC\ra \CC[G^{e,w_0}]$,
see~\cite[Theorem 2.10]{BeFoZe05}.

Note that by varying the reduced expression $\ii$ we obtain a priori several 
cluster algebra structures on $\CC[G^{e,w_0}]$, but according to
\cite[Remark 2.14]{BeFoZe05} all these structures coincide and give
rise to the same cluster variables and clusters. 
Note also that in type $\Aa_n$, the upper cluster algebra
$\overline{\mathcal{A}}(\ii)_\CC$ coincides with the cluster algebra
$\mathcal{A}(\ii)_\CC$, see~\cite[Remark 2.18]{BeFoZe05}.

\subsection{} \label{ssec:clustN}
Let $\tilde{\mathbf x}'$ be the subset 
of $\tilde{\mathbf x}$ obtained by removing 
the variables indexed by the $n$ non-exchan\-gea\-ble elements in $[1,r]$. 
Let $\mathcal{F}'$ be the field of rational functions over 
$\CC$ in the $r$ variables of $\tilde{\mathbf{x}}'$.
Finally, let $\widetilde{B}(\ii)'$ be
the matrix obtained from $\widetilde{B}(\ii)$ by removing 
the rows labelled by the $n$ non-exchangeable elements in $[1,r]$,
and let $\overline{\mathcal{A}}(\ii)'_\CC$ denote the upper cluster algebra
associated to the seed $(\tilde{\mathbf{x}}',\widetilde{B}(\ii)')$,
a subalgebra of $\mathcal{F'}$.
By \ref{1;2}, we see that the isomorphism of fields 
$\varphi'_\ii : \mathcal{F}' \ra \CC(N)$ defined by
\[
\varphi_\ii'(x_k)=\Delta(k,\ii), \qquad (x_k\in \tilde{\mathbf x}'),
\] 
restricts to an algebra isomorphism 
$\overline{\mathcal{A}}(\ii)'_\CC\ra \CC[N]$.

Clearly the same remarks as in \ref{1;3;4} apply to the cluster
algebra structures of $\CC[N]$. 

\subsection{}
We shall now describe in representation theoretic terms 
the restriction to $N$ of the regular function
$\Delta(k,\ii)$.
Let $L(\vpi_i)$ be the fundamental irreducible $\g$-module with 
highest weight $\vpi_i$. 
Fix a highest weight vector $u_{\vpi_i}$.

It is well-known that $L(\vpi_i)$ can be realized in a canonical
way as a subspace of the vector space $\CC[N]$ by restricting the summand
$L(\vpi_i)$ of $\CC[N_-\bcsl G]$ to $\CC[N]$. 
Thus, the highest weight
vector $u_{\vpi_i}$ becomes identified with the constant regular
function $\mathbf{1}$.
Using this identification we obtain the following lemma:

\begin{lemu}\label{lem11}
For $k\in [-n,-1] \cup [1,r]$ we have 
\[
\Delta(k,\ii) = \begin{cases}
f_{i_r}^{\max}f_{i_{r-1}}^{\max}\cdots f_{i_{k+1}}^{\max}(u_{\vpi_{i_k}}) 
&\text{ if } k\in [1,r],\\
f_{i_r}^{\max}f_{i_{r-1}}^{\max}\cdots f_{i_{1}}^{\max}(u_{\vpi_{\abs{i_k}}})
&\text{ if } k\in [-n,-1].
\end{cases}
\]
In particular for $-k\in [-n,-1]$, the minor 
$\Delta(-k,\ii) = \Delta_{\vpi_k,w_0(\vpi_k)}$
is a lowest weight vector of $L(\vpi_k)$.
\end{lemu}

\begin{pf} This follows from~\cite[p. 150--151]{BerZel97} 
and~\ref{lem:well-known} by restricting regular
functions on $G$ to the subgroup $N$, see also~\cite[p. 113]{BerZel01}.
\end{pf}

\section{Cluster variables and semicanonical basis} \label{sec:cv-scb}
\subsection{}
Let $Q$ be a Dynkin quiver such that $\abs{Q}$ is the diagram of $G$.
In this section we prove that if $\ii$ is a reduced expression for $w_0$
which is adapted to $Q$, 
the minors $\Delta(k,\ii)\in\CC[N]$ coincide with certain
dual semicanonical basis vectors coming from the injective modules
of the  Auslander algebra of $\CC Q$.  
In particular, this shows that the set of minors 
$\{\Delta(k,\ii)\}$ depends only on $Q$, not on the choice of
a particular expression $\ii$ adapted to $Q$.  

\subsection{}
Recall that by pushing down the injective modules
of the  Auslander algebra of $\CC Q$ we obtain a set 
of indecomposable rigid modules $\sM(j,\ii)$ over the preprojective
algebra $\Lambda$, see~\ref{ssec:startM} and~\ref{intro:expr}.
Since these modules are rigid, they have an open orbit in their
module variety, and the closure of this orbit is an irreducible
component. Therefore the module $\sM(j,\ii)$ can be used to label an element 
$\rho_{\sM(j,\ii)}$ of the dual semicanonical basis (see \cite{Lusz00},
\cite[Section 7.2]{GeLeSc04a}).

\subsection{}
The proof of Theorem~\ref{th1} will make use of certain results
of~\cite{GeLeSc04b} that we shall now recall.
Let $I_i$ denote the injective envelope of the simple $\Lambda$-module
$S_i$ with dimension vector $\de_i$, thus $I_i=D\Lam(-,i)$.
In~\cite{GeLeSc04b} the fundamental $\g$-module $L(\vpi_i)$ was
realized in terms of the lattice of submodules of $I_i$.
This goes as follows.

Let $\M$ denote Lusztig's algebra of constructible functions
on the varieties of finite-dimensional $\Lambda$-modules, 
and let $\M^*$ be its graded Hopf dual, an algebra isomorphic
to $\CC[N]$.
For a finite-dimensional $\Lambda$-module $X$, let 
$\delta_X$ denote the linear form on $\M$ 
obtained by evaluation at~$X$. Then in the identification
$\CC[N] \cong \M^*$, the subspace $L(\vpi_i)$ gets identified
to the subspace of $\M^*$ spanned by the linear forms $\delta_X$
where $X$ runs over the lattice of submodules of $I_i$,
and one has explicit formulas for the action of the Chevalley
generators of $\g$ on each vector $\delta_X$ \cite[Theorem~3]{GeLeSc04b}.
In particular $\delta_{I_i}=\rho_{I_i}$ is a lowest weight
vector of $L(\vpi_i)$, and $\delta_{0}$, where $0$ means the
zero submodule of $I_i$, is a highest weight vector. 

Let $X$ be a submodule of $I_i$. We have a short exact sequence
of $\Lambda$-modules 
\[
0 \ra X \ra I_i \xrightarrow{p} Y \ra 0,
\] 
where $Y$ is  determined  up to isomorphism by the isomorphism class 
of $X$, see~\cite[Lemma 1]{GeLeSc04b}. 
For $j\in [1,n]$ let $m_j$ denote the multiplicity of
$S_j$ in the socle of $Y$, and let $X_j$ be the unique submodule
of $I_i$ such that $X\subset X_j\subset I_i$ and
$X_j/ X$ is isomorphic to $S_j^{\oplus m_j}$.  Thus $X_j$ is the
pullback of $p$ and the inclusion of $S_j^{\oplus m_j}$ into $Y$.

\begin{lem}\label{lem32}
With the above notation, we have 
$f_j^{\max} (\delta_X) = f_j^{(m_j)}(\delta_X) = \delta_{X_j}$.  
\end{lem}
\begin{pf}
By~\cite[Theorem 3 (ii)]{GeLeSc04b},
we have that 
\[
f_j^k (\delta_X) = \int_{\f=(X=X(1) \subset \cdots \subset X(k))} 
\delta_{X(k)}
\]
where the integral is over the variety of flags $\f$ of submodules of $I_i$
such that $X(s)/ X(s-1)$ is isomorphic to $S_j$ for all $1<s\leq k$.
For $k>m_j$ this variety is empty by definition of $m_j$,
hence $f_j^k(\delta_X) = 0$.
For $k=m_j$, all flags $\f$ have their last step equal to 
$X_j$. Moreover since $X_j/X \cong S_j^{\oplus m_j}$
this variety is isomorphic to the variety of complete flags
in $\CC^{m_j}$, whose Euler characteristic is $m_j!$.
Hence $f_j^{m_j} (\delta_X) = m_j!  \delta_{X_j}$,
as claimed.
\end{pf}

\begin{lem}\label{lem33}
Let $k\in [-n,-1] \cup e(\ii)$ and $j=\theta^{-1}(k)$.
Then, if $i=\abs{i_k}$, the module $\sM(j,\ii)$ is a submodule
of $I_i$ and in $L(\vpi_i)$ there holds
\[
\delta_{\sM(j,\ii)}=\begin{cases}
f_{i_r}^{\max}f_{i_{r-1}}^{\max}\cdots f_{i_{k+1}}^{\max}(\delta_0)
&\text{ if } k\in e(\ii),\\
f_{i_r}^{\max}f_{i_{r-1}}^{\max}\cdots f_{i_{1}}^{\max}(\delta_0)
&\text{ if } k\in [-n,-1].
\end{cases}
\]
\end{lem}
\begin{pf}
Recall that  we defined the function $\theta$ in~\ref{LieG}.
If $k=-i\in [-n,-1]$, then  $\sM(j,\ii)=I_i$.
On the other hand in this case 
the product $f_{i_r}^{\max}f_{i_{r-1}}^{\max} \cdots 
f_{i_1}^{\max}$ maps $\delta_0$ to the lowest weight
vector of $L(\vpi_i)$, that is, to $\delta_{I_i}$,
as required.

If $k\in e(\ii)$, then $\RA_Q(x(j))$ belongs to the $\tau$-orbit
of $\RA_Q(x(k))$ by definition of $\theta$, therefore to the
$\tau$-orbit of $D\kQ(-,i)$ see~\ref{intro:expr}. 
It follows  that $\sM(j,\ii)$ is a 
submodule of $I_i$,
see~\ref{ssec:startM}.

More precisely, $j=\theta^{-1}(k)=\min\{l\in [k+1,r]\mid i_l=i\}$. We conclude
that in $L(\vpi_i)$ holds 
$f_{i_{j-1}}^{\max}\cdots f_{i_{k+1}}^{\max}(\delta _0)=\delta_0$ 
by~\ref{lem32}, since the socle of $I_i$ is $S_i$. 

Now consider for $l\in [j-1,r]$ the $\tLam$-submodule
$\widetilde{I}_{x(j)}(\leq l,\ii)$ of $D\tLam(-,x(j))$ with
\[
\widetilde{I}_{x(j)}(\leq l,\ii)(x(m)):=\begin{cases}
D\tLam(x(m),x(j)) &\text{ if } m\leq l,\\ 0 &\text{ else.}
\end{cases}
\]
With $I_{x(j)}(\leq l,\ii):=F_\lam \widetilde{I}_{x(j)}$ we see that
$I_{x(j)}(\leq r,\ii)=\sM(j,\ii)$ is a 
submodule of $I_i$ since $F_\lam D\tLam(-,x(j))\cong I_i$.

Using Lemma~\ref{lem32} we conclude that in $L(\vpi_i)$ we have
\[
f_{i_l}^{\max}(\delta_{I_{x(j)}(\leq l-1,\ii)})=\delta_{I_{x(j)}(\leq l,\ii)}
\text{ for } l\in[j,r].
\]
This yields that $\delta_{\sM(j,\ii)}$
is the extremal vector of $L(\vpi_i)$ with weight
$s_{i_r}s_{i_{r-1}}\cdots s_{i_{k+1}}(\vpi_i)$,
and the lemma follows.
\end{pf}

\bigskip
We can now finish the proof of Theorem~\ref{th1}.
Using Lemma~\ref{lem11} and Lemma~\ref{lem33} we obtain that
$\Delta(k,\ii)=\delta_{\sM(j,\ii)}$.
But since $\sM(j,\ii)$ is rigid, its orbit is open and we
have $\delta_{\sM(j,\ii)} = \rho_{\sM(j,\ii)}$.

\section{Examples}
Our running example will be the Dynkin quiver $Q$ of type $\Dd_5$ with
the following orientation:
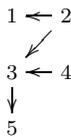
\[
\def\objectstyle{\scriptstyle}
\xymatrix@-1.2pc @l{
&5\\
4\ar[r]&3\ar[u]\\
2\ar[ru]\ar[r]& 1}
\]

\subsection{$\AuC{Q}$ for $\Dd_5$.}
\label{expl-GaD5}
We show the quiver $A_Q$ of $\AuC{Q}$:
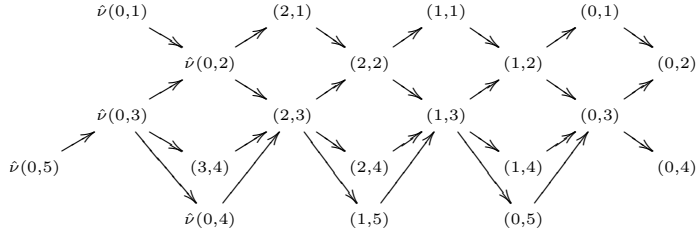
\[
\def\objectstyle{\scriptscriptstyle}
\xymatrix@-1.4pc @l{
        &     &(0,5)\ar[ldd]&  &(1,5)\ar[ldd]&    &\hnu(0,4)\ar[ldd]\\
(0,4)   &     &(1,4)\ar[ld] &  &(2,4)\ar[ld] &&(3,4)\ar[ld]&&\ar[ld]\hnu(0,5)\\
 &(0,3)\ar[ld]\ar[lu]&&(1,3)\ar[lu]\ar[luu]\ar[ld]&&(2,3)\ar[lu]\ar[luu]\ar[ld]&&\hnu(0,3)\ar[lu]\ar[ld]\ar[luu]\\
(0,2) &&(1,2)\ar[lu]\ar[ld]&&(2,2)\ar[lu]\ar[ld]&&\hnu(0,2)\ar[lu]\ar[ld]\\
&(0,1)\ar[lu]&&(1,1)\ar[lu]&&(2,1)\ar[lu]&&\hnu(0,1)\ar[lu]
}
\]
The isomorphism classes of the indecomposable representations of
$Q$ correspond to the vertices of $A_Q$. Each such representation
is uniquely determined by its dimension vector:
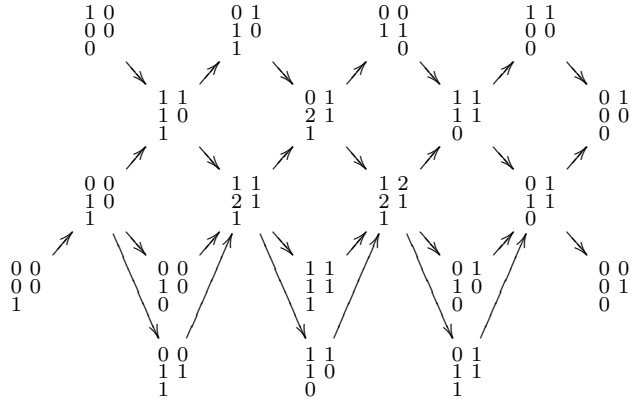
\[
\def\objectstyle{\scriptscriptstyle}
\xymatrix@-1.4pc @l{
&&{\bsm 0&1\\1&1\\1& \esm}\ar[ldd]&&{\bsm 1&1\\1&0\\ 0\esm}\ar[ldd]&&{\bsm 0&0\\1&1\\ 1\esm}\ar[ldd]\\
{\bsm 0&0\\0&1\\ 0\esm}&&{\bsm 0&1\\1&0\\ 0\esm}\ar[ld]&&{\bsm 1&1\\1&1\\ 1\esm}\ar[ld]&&{\bsm 0&0\\1&0\\ 0\esm}\ar[ld]&&{\bsm 0&0\\0&0\\ 1\esm}\ar[ld]\\
 &{\bsm 0&1\\1&1\\ 0\esm}\ar[lu]\ar[ld]&&{\bsm 1&2\\2&1\\ 1\esm}\ar[lu]\ar[luu]\ar[ld]&&{\bsm 1&1\\2&1\\ 1\esm}\ar[lu]\ar[luu]\ar[ld]&&{\bsm 0&0\\1&0\\ 1\esm}\ar[lu]\ar[ld]\ar[luu]\\
{\bsm 0&1\\0&0\\ 0\esm}&&{\bsm 1&1\\1&1\\ 0\esm}\ar[lu]\ar[ld]&&{\bsm 0&1\\2&1\\ 1\esm}\ar[lu]\ar[ld]&&{\bsm 1&1\\1& 0\\ 1\esm}\ar[lu]\ar[ld]\\
&{\bsm 1&1\\0&0\\ 0\esm}\ar[lu]&&{\bsm 0&0\\1&1\\&0\esm}\ar[lu]&&{\bsm 0&1\\1&0\\ 1\esm}\ar[lu]&&{\bsm 1&0\\0&0\\ 0\esm}\ar[lu]
}
\]
The dimension vector of the injective $\AuC{Q}$-module $D\AuC{Q}(-,(0,q))$
is given by the $q$-com\-po\-nents of the corresponding dimension vectors. 
The dimension vector of $D\AuC{Q}(-,\tau^j(0,q))$ is obtained from this
by ``translation and cut-off''.
Here we show for example $\dimv D\AuC{Q}(-,(0,3))$ and
$\dimv D\AuC{Q}(-,\tau(0,3))$:
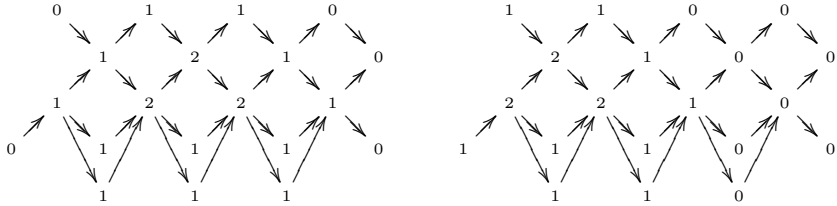
\[
\def\objectstyle{\scriptscriptstyle}
\xymatrix@-1.4pc @l{
         &     &1\ar[ldd]&      &1\ar[ldd]&       &1\ar[ldd]\\
0        &     &1\ar[ld] &      &1\ar[ld] &       &1\ar[ld]&&0\ar[ld]\\
 &1\ar[lu]\ar[ld]&&2\ar[lu]\ar[luu]\ar[ld]&&2\ar[lu]\ar[luu]\ar[ld]&&1\ar[lu]\ar[luu]\ar[ld]\\
0&&1\ar[lu]\ar[ld]&&2\ar[lu]\ar[ld]&&1\ar[lu]\ar[ld]\\
&0\ar[lu]&&1\ar[lu]&&1\ar[lu]&&0\ar[lu]
}
\qquad
\def\objectstyle{\scriptscriptstyle}
\xymatrix@-1.4pc @l{
    &&\ar[ldd]0&  &\ar[ldd]1&&\ar[ldd]1\\
0   &&\ar[ld]0 &  &\ar[ld] 1&&\ar[ld] 1&&\ar[ld]1\\
&\ar[lu]\ar[ld]0&&\ar[lu]\ar[luu]\ar[ld]1&&\ar[luu]\ar[lu]\ar[ld]2&&\ar[lu]\ar[luu]\ar[ld]2\\
0&&\ar[lu]\ar[ld]0&&\ar[lu]\ar[ld]1&&\ar[lu]\ar[ld]2\\
&\ar[lu]0&&\ar[lu]0&&\ar[lu]1&&\ar[lu]1
}
\]

\subsection{The category $\cGam_Q$.}
\label{expl-ga*}
We display here the quiver of the (graded) category $\cGam_Q$ associated
to a quiver of type $\Dd_5$ with the same orientation as above. The arrows
of degree one are dotted.
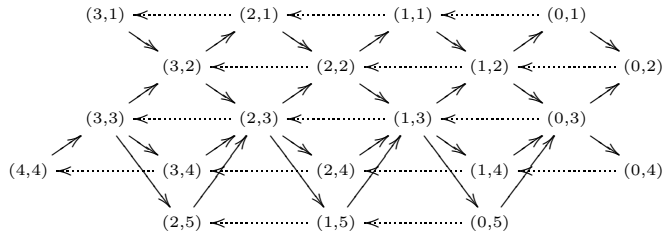
\[
\def\objectstyle{\scriptscriptstyle}
\xymatrix@-1.4pc @l{
&&\ar[ldd](0,5)\ar@{.>}[rr]&&\ar[ldd](1,5)\ar@{.>}[rr]&&\ar[ldd](2,5)\\
(0,4)\ar@{.>}[rr]&&\ar[ld](1,4)\ar@{.>}[rr]&&\ar[ld](2,4)\ar@{.>}[rr]&&\ar[ld](3,4)\ar@{.>}[rr]&&\ar[ld](4,4)\\
&\ar[lu]\ar[ld](0,3)\ar@{.>}[rr]&&\ar[lu]\ar[luu]\ar[ld](1,3)\ar@{.>}[rr]&&\ar[lu]\ar[luu]\ar[ld](2,3)\ar@{.>}[rr]&&\ar[lu]\ar[luu]\ar[ld](3,3)\\
(0,2)\ar@{.>}[rr]&&\ar[lu]\ar[ld](1,2)\ar@{.>}[rr]&&\ar[lu]\ar[ld](2,2)\ar@{.>}[rr]&&\ar[lu]\ar[ld](3,2)\\
&\ar[lu](0,1)\ar@{.>}[rr]&&\ar[lu](1,1)\ar@{.>}[rr]&&\ar[lu](2,1)\ar@{.>}[rr]&&\ar[lu](3,1)
}
\]
Recall that also the relations are easy to read off: For each arrow
$\alpha\df x\ra y$ we have the corresponding mesh relation~\eqref{eqn:mesh-rel}
from $y$ to $x$ if $\alpha$ is of degree $1$. Otherwise,
if  $x$  is not an ``injective'' vertex, 
(i.e.~here if $x\not\in\{(0,1),(0,3),(0,5)\}$) there are
one or two paths of length $2$ (and degree $1$) from $y$ to $x$. 
The first case occurs when $y$ is a ``projective'' vertex 
(i.e.~here if $y\in\{(3,2),(3,3),(2,5)\}$) and
the unique path of length $2$ from $y$ to $x$ is a zero-relation, otherwise
the two paths form a commutativity relation. 

\subsection{A projective-injective $\cGam_Q$-module.}
\label{expl-ga*-pi}
We display for $\cGam_Q$ as in~\ref{expl-ga*} the projective
module $\cGam_Q((0,5),-)$. Since $(0,5)$ is an injective vertex
of $\AR{Q}$ we have $\cGam_Q((0,5),-)\cong D\cGam_Q(-,(0,\mu(5))$. Moreover,
the projective modules $\cGam_Q((j,5),-)$ for $0\leq j\leq 2$ are
easily found as submodules of  $\cGam_Q((0,5),-)$.

Each entry $(n,q)$ represents 
a basis vector which corresponds to a (graded) simple composition factor of 
this type. The arrows indicate as usual the action of $\cGam_Q$.
\pagebreak[2]
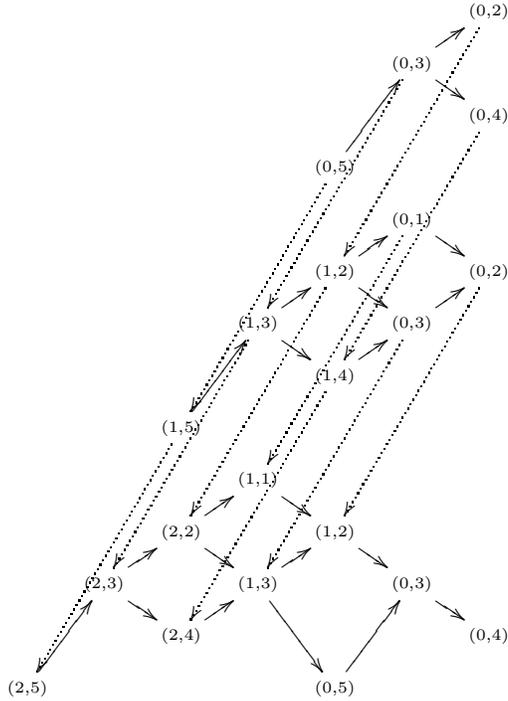
\[
\def\objectstyle{\scriptscriptstyle}\def\labelstyle{\scriptscriptstyle}
\xymatrix@-1.4pc @l{
&&\ar[ldd](0,5)&&&&\ar[ldd](2,5)\\
(0,4)&&&&\ar[ld](2,4)\\
&\ar[lu](0,3)&&\ar[luu]\ar[ld](1,3)&&\ar[lu]\ar[ld](2,3)\\
&&\ar[lu](1,2)&&\ar[lu]\ar[ld](2,2)\\
&&&\ar[lu](1,1)\\
&&&&\ar[ldd](1,5)\ar@{.>}[rruuuuu]\\
&&\ar[ld](1,4)\ar@{.>}[rruuuuu]\\
&\ar[ld](0,3)\ar@{.>}[rruuuuu]&&\ar[lu]\ar[ld](1,3)\ar@{.>}[rruuuuu]\\
(0,2)\ar@{.>}[rruuuuu]&&\ar[lu]\ar[ld](1,2)\ar@{.>}[rruuuuu]\\
&\ar[lu](0,1)\ar@{.>}[rruuuuu]\\
&&\ar[ldd](0,5)\ar@{.>}[rruuuuu]\\
(0,4)\ar@{.>}[rruuuuu]\\
&\ar[lu]\ar[ld](0,3)\ar@{.>}[rruuuuu]\\
(0,2)\ar@{.>}[rruuuuu]\\
}
\]

\subsection{Dimensions.}
We include the result of some calculations of 
\[
d(Q):=\dim\End_\Lam(\sM_Q)
\]
for specific orientations. This can be done quite easily on a computer
using the formula~\eqref{eq:endim}.
\begin{alignat*}{2}
A_n &\df
{\def\objectstyle{\scriptstyle}\vcenter{
\xymatrix @-1.2pc{1\ar[r]&2\ar[r]&\quad\cdots\quad\ar[r]&(n-1)\ar[r]&n}}}
&d(A_n)&= 2 \binom{n}{5}+7\binom{n}{4}+9\binom{n}{3}+5\binom{n}{2}+n,\\
D_n&\df
\def\objectstyle{\scriptstyle}
\vcenter{\xymatrix @-1.2pc{
       &         &                      &                   &(n-1)\\
1\ar[r]& 2\ar[r] &\quad\cdots\quad\ar[r]&(n-2)\ar[ru]\ar[rd]\\
       &         &                      &                   & n}}\quad
&
d(D_n) &= 27 \binom{n}{5} + 43\binom{n}{4}+ 19\binom{n}{3}+ 2\binom{n}{2},\\
 E_n&\df
 \def\objectstyle{\scriptstyle}
 \vcenter{\xymatrix@-1.2pc{
         &           &                      &2 \ar[rd]\\
         &           &               4\ar[r]&3 \ar[r] & 1\\
 n\ar[r] &(n-1)\ar[r]&\quad\cdots\quad\ar[r]&5\ar[ru]}}
&
d(E_n)&=\begin{cases}   
  2444 &\text{ if } n=6,\\
 13130 &\text{ if } n=7,\\
107114 &\text{ if } n=8.
\end{cases}
\end{alignat*}

\subsection{An adapted ordering on $\Obj(\AuC{Q})$:} \label{expl:ada-ord}
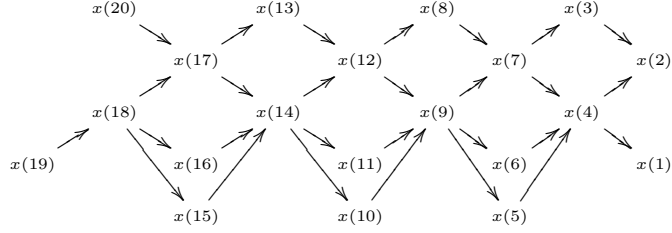
\[
\def\objectstyle{\scriptscriptstyle}
\xymatrix@-1.4pc @l{
&&x(5)\ar[ldd]&&x(10)\ar[ldd]&&x(15)\ar[ldd]\\
x(1)&&\ar[ld]x(6)&&\ar[ld]x(11) &&\ar[ld]x(16)&&x(19)\ar[ld]\\
 &\ar[lu]\ar[ld]x(4)&&\ar[lu]\ar[luu]\ar[ld]x(9)&&\ar[lu]\ar[luu]\ar[ld]x(14)&&\ar[lu]\ar[luu]\ar[ld]x(18)\\
x(2)&&\ar[lu]\ar[ld]x(7)&&\ar[lu]\ar[ld]x(12)&&\ar[lu]\ar[ld]x(17)\\
&x(3)\ar[lu]&&\ar[lu]x(8)&&\ar[lu]x(13)&&\ar[lu]x(20)
}
\]
According to~\ref{intro:expr} we obtain for $w_0$ the following 
reduced expression which is adapted to $Q$:
\[
\ii=(4,2,1,3,5,4,2,1,3,5,4,2,1,3,5,4,2,3,4,1)
\]
\subsection{The quiver $\tAi$.} \label{expl:Gam-i}
For the adapted expression $\ii$ from~\ref{expl:ada-ord} we obtain
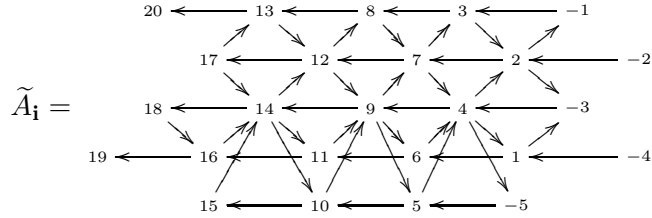
\[
\def\objectstyle{\scriptscriptstyle}
\tAi=\vcenter{\xymatrix@-1.4pc @l{
&&{-5}\ar[rr]&&{5}\ar[rr]\ar[ldd]&&{10}\ar[rr]\ar[ldd]&&{15}\ar[ldd]\\
{-4}\ar[rr]&&\ar[ld]{1}\ar[rr]&&\ar[ld]{6}\ar[rr]&&\ar[ld]{11}\ar[rr]&&\ar[ld]{16}\ar[rr]&&{19}\\
 &{-3}\ar[rr]&&\ar[lu]\ar[ld]\ar[luu]{4}\ar[rr]&&\ar[lu]\ar[luu]\ar[ld]{9}\ar[rr]&&\ar[lu]\ar[luu]\ar[ld]{14}\ar[rr]&&\ar[lu]{18}\\
{-2}\ar[rr]&&\ar[lu]\ar[ld]{2}\ar[rr]&&\ar[lu]\ar[ld]{7}\ar[rr]&&\ar[lu]\ar[ld]{12}\ar[rr]&&\ar[lu]\ar[ld]{17}\\
&{-1}\ar[rr]&&{3}\ar[lu]\ar[rr]&&\ar[lu]{8}\ar[rr]&&\ar[lu]{13}\ar[rr]&&{20}
}}
\]
The exchangeable vertices are $\{1,2,3,4,5,6,7,8,9,10,11,12,13,14,16\}$

\section{Appendix: Using the triangulated structure} 
\label{sec:apx}
\begin{prp}[Freyd/Heller] 
\label{prp:FK}
Let $\cD$ be a triangulated $\ka$-category with suspension functor $\Sig$.
\begin{itemize}
\item[(a)]
The category $\lmd{\cD}$ is a Frobenius category. Thus the stable category
$\slmd{\cD}$ is triangulated with suspension functor $\Ome_\cD^{-1}$, 
the inverse of Heller's loop functor.
\item[(b)]
In $\slmd{\cD}$ we have a functorial isomorphism
$M^{\Sig}\cong \Ome_\cD^{3} M$.
\end{itemize}
\end{prp}

Part (a) is from~\cite[Section 3]{Freyd66a}, see 
also~\cite[Chapter 5]{Neeman01} for a modern treatment.
Part (b) is a special case of~\cite[\S 16]{Heller68a}, see also
\cite[Proposition B.2]{Krau02}.

\begin{rems} \label{ssec:cov-rem1}
(a) We may consider $\cD$ as a $\cD\text{-}\cD$-bimodule,
i.e. a functor $\cD^{\text{op}}\times\cD\ra\lmd{\ka}$.
Similarly, $D\cD$ with $D\cD(a,b):=\Hom_\ka(\cD(b,a),\ka)$ is a 
$\cD\text{-}\cD$-bimodule.

(b) Suppose that $\cD$ admits Auslander-Reiten triangles with translate 
$\tau\df\cD\ra\cD$. In this case we set $\nu:=\Sig\tau$. Then
the Auslander-Reiten formula $\cD(x,\Sig y)\cong \Hom_k(\cD(y,\tau x),\ka)$
may be interpreted as an isomorphism of
bimodules
\begin{equation} \label{eqn:cov-AR}
D\cD \cong \cD^{\nu^{-1}} \quad\text{where } 
\cD^{\nu^{-1}}(a,b):=\cD(a,\nu b).
\end{equation}
We conclude that 
\[
\Nak M:= D\cD\otimes_\cD M\cong \cD^{\nu^{-1}}\otimes_\cD M\cong M^{\nu^{-1}}
\]
is a Nakayama functor for $\lmd{\cD}$ and $\nu^{-1}$ the corresponding
Nakayama automorphism for $\cD$, see for example~\cite[Section~2]{gabr80}.
In this situation we will write
\[
M^{(i)}:= M^{\tau^i}.
\]

(c) If moreover $\cD$ is locally bounded, then $\slmd{\cD}$  admits 
Auslander-Reiten triangles with translation $\tau_\cD=\Ome_\cD^2\Nak$.
With Proposition~\ref{prp:FK} we obtain the  functorial isomorphisms
\begin{equation} \label{eqn:cov-form}
M^{(-1)}\cong \tau_\cD\Ome_\cD M\quad\text{ and }\quad 
\tau_\cD^{-3}M\cong M^{\Sig\tau^3} \text{ (in } \slmd{\cD}\text{)}.
\end{equation}
In fact, we have 
$\tau_\cD\Ome_D\cong\Ome^3\Nak\cong \,?^{\Sig\nu^{-1}}=\,?^{\tau^{-1}}$ and
$\tau_\cD^{-3}=\Ome_\cD^6\Nak^3\cong\,?^{\Sig^{-2}\nu^3}\cong\,?^{\Sig \tau^3}$.
(d) If we  consider in this context {\em right}  modules 
(i.e.~contravariant functors),
we get in $\srmd{\cD}$ a functorial isomorphism $M^{\Sig}\cong\Ome_\cD^{-3} M$
and consequently
\begin{equation*} 
M^{(1)}\cong \tau_\cD\Ome_\cD M\quad\text{ and } 
\tau_\cD^{3}M\cong M^{\Sig\tau^3} \text{ (in } \srmd{\cD}\text{)}.
\end{equation*}
\end{rems}

\subsection{Derived categories.} \label{apx:derc}
It follows from Happel's 
description~\cite[I.5.6]{happ88} of the derived category 
$\Der{\kQ\op}:=\Der{\lmd{\kQ\op}}$ that we have a natural equivalence 
\[
\linj{\tLam}\cong  \Der{\kQ\op}.
\]
In particular, $\linj{\tLam}$ is a triangulated category which admits
Auslander-Reiten triangles.
The suspension functor resp. the Auslander-Reiten translate are
\[
\Sig I= I^{\hnu\tau}\quad\text{ resp. } \tau I=I^{\tau^{-1}},
\]
see~\ref{ssec:a-conv}.
In our situation, these functors are  up to isomorphism determined by 
their effect on objects.
We conclude that we have an isomorphism of bimodules
\begin{equation} \label{eqn:bi-DtL}
\tLam^{\hnu}\cong D\tLam,
\end{equation}
see~\ref{ssec:cov-rem1}~(b).

In order to state our next result we introduce $\uGam_Q$, the full subcategory
of the Auslander category $\AuC{Q}$ which contains all objects except those of 
the form  $\hnu (0,q)$ for $q\in Q_0$. 
We call $\uGam_Q$ the {\em stable Auslander category} of $\kQ$.

\begin{prp} \label{apx:repc}
The category $\tLam$ is isomorphic to the repetitive category 
$\huGam_Q$ of the stable Auslander category of $\kQ$.
\end{prp}

\begin{pf} Recall that by definition the objects of $\huGam_Q$ are of the 
form $(z,x)$ with $z\in\Zet$ and $x\in\Obj(\uGam_Q)$ and we have
\[
\huGam_Q((z,x),(z',x'))=\begin{cases}
\uGam_Q(x,x')  &\text{ if } z=z',\\
D\uGam_Q(x,x') &\text{ if } z=z'-1,\\
0            &\text{ else.}
\end{cases}
\]
Here we define the dual $\uGam_Q\text{-}\uGam_Q$-bimodule 
$D\uGam_Q$ by $D\uGam_Q(x,y):=\Hom_\ka(\uGam_Q(y,x),k)$, 
compare~\ref{ssec:cov-rem1}~(a). Now, by the bimodule 
isomorphism~\eqref{eqn:bi-DtL} we see that
the assignation
\[
(z, (i,q))\mapsto \tau^i\hnu^{-z}(0,q)
\]
induces an isomorphism $\huGam_Q\ra\tLam$.
\end{pf}

\subsection{Conclusions.}
(a) In our situation we note that in $\Der{\kQ}\cong\Der{\kQ\op}$ we have 
an isomorphism of functors 
\begin{equation} \label{eqn:susptau}
\Sig^2\cong \tau^{-h(Q)}
\end{equation}
where $h(Q)$ is the Coxeter number of $\abs{Q}$. 
In fact, both functors coincide on
objects as one easily verifies in $\tLam$. In our quiver situation this is 
sufficient. From~\eqref{eqn:cov-form} we obtain immediately  the remarkable
functorial isomorphism
\begin{equation}\label{eqn:cov-tau}
\tau_{\tLam}^6 \cong M^{(h(Q)-6)}
\end{equation}
in $\slmd{\tLam}\cong\slmd{\Der{\kQ}}$. 

(b) The action of the infinite cyclic group $\ebrace{\tau}$ on $\tLam$ provides
us with a {\em Galois covering}
\[
F\df \tLam\ra \Lam,
\]
see for example~\cite{gabr81}.
We conclude that $\linj{\Lam}$ is the orbit category of $\linj{\tLam}$ modulo
the induced action of $\ebrace{\tau}$. Now, $\linj{\tLam}\cong\Der{\kQ}$
is a triangulated category, and the hypothesis of~\cite{Kell05} are obviously
fulfilled. Thus $\linj{\Lam}$  is also a triangulated category with
Auslander-Reiten triangles and the corresponding translation $\bar{\tau}$
is the identity. Note moreover that the induced suspension is isomorphic
to the corresponding Nakayama automorphism, i.e. $\bar{\Sig}\cong\bar{\nu}$. 
By applying~\eqref{eqn:cov-form} to $\cD=(\linj{\Lam})\op$ we conclude
that 
\[
\tau_\Lam^3 M\cong M^{\bar{\nu}}\quad\text{and}\quad\tau_\Lam^{6}M\cong M
\]
holds functorially in $\slmd{\Lam}$. In case $\nu$ is just a translation,
see~\ref{ssec:cov-constr}, we even have  $\tau^3_\Lam M\cong M$.
This is our interpretation of the proof for the $6$-periodicity of 
$\tau_\Lam$ in~\cite{AusRei96}. Even more directly by~\eqref{eqn:cov-form}
we conclude that
\[
\tau_\Lam\Ome_\Lam M\cong M
\]
functorially. This means that the triangulated category $\slmd{\Lam}$ is
of Calabi-Yau dimension~2. 

(c) Since $\tLam\cong \huGam_Q$ we have by Happel's 
Theorem~\cite[II.4]{happ88} $\slmd{\tLam}\cong \Der{\lmd{\uGam_Q}}$ as
triangulated categories. If we consider the push-down functor
\[
F_\lam\df \slmd{\tLam}\ra \slmd{\Lam}
\]
associated to the Galois covering $F\df\tLam\ra \Lam$ 
the subcategory of  $\Lam$-modules of the first kind 
(i.e. the subcategory of objects which are isomorphic to a push-down) 
is  equivalent to the 
orbit category $\Der{\lmd{\uGam_Q}}/ \ebrace{\tau\Sig^{-1}}$ via the above 
identifications and~\eqref{eqn:cov-form}, see~\cite{dosk87}. Here,
$\Sig$ resp.~$\tau$  are the suspension resp.~the Auslander-Reiten translation
in $\Der{\lmd{\uGam_Q}}$.

Now, for a Dynkin quiver $Q$ with Coxeter number $h(Q)\leq 6$ the
algebras $\uGam_Q$ are (quasi-) tilted. Thus, in these cases we find 
$\slmd{\Lam}\cong\Der{\uGam_Q}/\ebrace{\Sig\tau^{-1}}$ is a cluster
category in the sense of~\cite{BMRRT04}.

\subsection{Remark} Let $H$ be a (basic, connected) finite-dimensional
hereditary $\ka$-algebra of finite representation type. Then $H$ is a
species of type $\Aa$ -- $\Gg$ in the sense of Dlab and Ringel, 
see~\cite{DlaRin76}. In this case we may  also study the stable Auslander
category $\uGam_H$. The same argument as in~\ref{apx:derc} and~\ref{apx:repc}
shows that $\linj{\huGam_H}\cong\Der{H^{\op}}$. The Auslander-Reiten
translate in $\Der{H^{\op}}$ induces an automorphism $\tau$ of
$\huGam_H$. Thus we may consider the Galois covering
$\huGam_H\ra \huGam_H/\ebrace{\tau}$. So we are tempted to consider 
$\huGam_H/\ebrace{\tau}$ as the preprojective algebra of $H$. However,
$\tau$ is now in general not determined by its effect on objects since
$\operatorname{Out}_\ka(H)=\operatorname{Aut}_k(H)/\operatorname{Inn}(H)$ 
is possibly
non-trivial. Thus we are (for the moment) unable to compare
$\huGam_H/\ebrace{\tau}$ with the possible choices for the preprojective
algebra of $H$ in the sense of Dlab and Ringel~\cite{DlaRin80}.

Anyway, if we denote by $\abs{H}$ the (unoriented) diagram of $H$ then
we find the list which we present in Figure~\ref{fig:types} 
of interest due to its similarity with the cluster types of $\CC[N]$.
In the case of Coxeter number $c(\abs{H})=6$ we display the diagram of a 
canonical tubular algebra following Lenzing~\cite{Lenz96} and
the corresponding extended affine root system in the sense of 
Saito~\cite{sait85}. In the case of $\Gg_2$ there are two different diagrams
of canonical algebras which  produce derived equivalent algebras
(for adequate choices of bimodules), 
the corresponding two root systems are isomorphic up to the marking.
\begin{figure} 
\caption{\label{fig:types}}
\vspace{2ex}
\begin{tabular}{c|c|l}
$c(\abs{H})$ & $\abs{H}$ & $\uGam_H$ is \\
             &           & tilted of\\\hline
   3         &  $\Aa_2$  &  $\Aa_1$ \\\hline
   4         &  $\Aa_3$  &  $\Aa_3$\\
             &  $\Bb_2$  &  $\Bb_2$\\\hline
   5         &  $\Aa_4$  &  $\Dd_6$\\\hline
\end{tabular}\qquad
\begin{tabular}{c|c|c|c}
$c(\abs{H})$ & $\abs{H}$ & $\uGam_H$ is tilted of& root system\\\hline
   6         &  $\Aa_5$  &  
$\vcenter{\UseComputerModernTips 
\xymatrix@-1.4pc{
&\cdot\ar[r]&\cdot\ar[r]&\cdot\ar[r]&\cdot\ar[r]&\cdot\ar[rd]\\
\cdot\ar[ru]\ar[rr]\ar[rrrd]&&\ar[rr]&&\cdot\ar[rr]&&\cdot\\
&&&\cdot\ar[rrru]}}$ & $\widetilde{\Ee}_8^{(1,1)}$\\
& $\Bb_3$ &
$\vcenter{\def\labelstyle{\scriptscriptstyle}\UseComputerModernTips
\xymatrix@-1.4pc{
&\cdot\ar[r]&\cdot\ar[rd]^{(2,1)}\\
\cdot\ar[ru]^{(2,1)}\ar[rd]&&&\cdot\\
&\cdot\ar[r]&\cdot\ar[ru]}}$  &$\widetilde{\Ff}_4^{(1,1)}$\\
& $\Cc_3$ &
$\vcenter{\def\labelstyle{\scriptscriptstyle}\UseComputerModernTips
\xymatrix@-1.4pc{
&\cdot\ar[r]&\cdot\ar[rd]^{(1,2)}\\
\cdot\ar[ru]^{(1,2)}\ar[rd]&&&\cdot\\
&\cdot\ar[r]&\cdot\ar[ru]}}$ & $\widetilde{\Ff}_4^{(2,2)}$\\
& $\Dd_4$ &
$\vcenter{\UseComputerModernTips
\xymatrix@-1.4pc{
&\cdot\ar[r]&\cdot\ar[rd]\\
\cdot\ar[ru]\ar[rd]\ar[r]&\cdot\ar[r]&\cdot\ar[r]&\cdot\\
&\cdot\ar[r]&\cdot\ar[ru]}}$  & $\widetilde{\Ee}_6^{(1,1)}$\\
& $\Gg_2$ &
$\vcenter{\def\labelstyle{\scriptscriptstyle}\UseComputerModernTips
\xymatrix@-1.4pc{
&\cdot\ar[r]&\cdot\ar[rd]^{(1,3)}\\
\cdot\ar[ru]^{(3,1)}&&&\cdot}}$ & $\widetilde{\Gg}_2^{(3,1)}$ \\
&  &
$\vcenter{\def\labelstyle{\scriptscriptstyle}\UseComputerModernTips
\xymatrix@-1.4pc{
&\cdot\ar[r]&\cdot\ar[rd]^{(3,1)}\\
\cdot\ar[ru]^{(1,3)}&&&\cdot}}$ & $\widetilde{\Gg}_2^{(1,3)}$ \\
\end{tabular}
\end{figure}
Note moreover, that in this case our previous calculations~\eqref{eqn:cov-tau}
predict that in the stable module category of the repetitive algebra
$\slmd{\huGam_H}$ the Auslander-Reiten translate should be $6$-periodic.
This is in fact the case for all (tubular) algebras which are derived 
equivalent to a canonical algebra with a diagram from our list. 
\subsection*{Acknowledgement} We like to thank Henning Krause for drawing
our attention to the result on functor categories over triangulated
categories. We also thank Markus Reineke for pointing out the reference for
the material in~\ref{intro:expr}.

\bibliography{ulitv}
\bibliographystyle{hamsplain} 
\end{document}